%% Plain TeX file
%%
%% Numerical characterization of the
%% K\"ahler cone of a compact K\"ahler manifold
%%
%% Jean-Pierre Demailly, Universit\'e Joseph Fourier, Grenoble
%% Mihai Paun, Universit\'e Louis Pasteur, Strasbourg

\magnification=1200
\hsize=12.5cm
\vsize=19.5cm
\parskip 3pt plus 1pt minus 1pt
\parindent=6.4mm

\let\sl\it

\font\eightrm=cmr8
\font\eightbf=cmbx8
\font\eightsy=cmsy8
\font\eighti=cmmi8
\font\eightit=cmti8
\font\eighttt=cmtt8
\font\eightsl=cmsl8

\font\sixrm=cmr6
\font\sixbf=cmbx6
\font\sixsy=cmsy6
\font\sixi=cmmi6

\font\fourteenbf=cmbx10 at 14.4pt
\font\twelvebf=cmbx10 at 12pt

\font\tenmsa=msam10
\font\sevenmsa=msam7
\font\fivemsa=msam5
\newfam\msafam
  \textfont\msafam=\tenmsa 
  \scriptfont\msafam=\sevenmsa
  \scriptscriptfont\msafam=\fivemsa

\font\tenmsb=msbm10
\font\eightmsb=msbm8
\font\sevenmsb=msbm7
\font\fivemsb=msbm5
\newfam\msbfam
  \textfont\msbfam=\tenmsb
  \scriptfont\msbfam=\sevenmsb
  \scriptscriptfont\msbfam=\fivemsb
\def\Bbb{\fam\msbfam\tenmsb}

\catcode`\@=11
\def\eightpoint{%
  \textfont0=\eightrm \scriptfont0=\sixrm \scriptscriptfont0=\fiverm
  \def\rm{\fam\z@\eightrm}%
  \textfont1=\eighti \scriptfont1=\sixi \scriptscriptfont1=\fivei
  \def\oldstyle{\fam\@ne\eighti}%
  \textfont2=\eightsy \scriptfont2=\sixsy \scriptscriptfont2=\fivesy
  \textfont\itfam=\eightit
  \def\it{\fam\itfam\eightit}%
  \textfont\slfam=\eightsl
  \def\sl{\fam\slfam\eightsl}%
  \textfont\bffam=\eightbf \scriptfont\bffam=\sixbf
  \scriptscriptfont\bffam=\fivebf
  \def\bf{\fam\bffam\eightbf}%
  \textfont\ttfam=\eighttt
  \def\tt{\fam\ttfam\eighttt}%
  \def\Bbb{\fam\msbfam\eightmsb}%
  \textfont\msbfam=\eightmsb
  \def\Cal{\fam\Calfam\eightCal}%
  \textfont\Calfam=\eightCal
  \abovedisplayskip=9pt plus 2pt minus 6pt
  \abovedisplayshortskip=0pt plus 2pt
  \belowdisplayskip=9pt plus 2pt minus 6pt
  \belowdisplayshortskip=5pt plus 2pt minus 3pt
  \smallskipamount=2pt plus 1pt minus 1pt
  \medskipamount=4pt plus 2pt minus 1pt
  \bigskipamount=9pt plus 3pt minus 3pt
  \normalbaselineskip=9pt
  \setbox\strutbox=\hbox{\vrule height7pt depth2pt width0pt}%
  \let\bigf@ntpc=\eightrm \let\smallf@ntpc=\sixrm
  \normalbaselines\rm}
\catcode`\@=12

\def\bC{{\Bbb C}}
\def\bN{{\Bbb N}}

\def\bQ{{\Bbb Q}}
\def\bR{{\Bbb R}}
\def\bZ{{\Bbb Z}}

\font\tenCal=eusm10
\font\eightCal=eusm8
\font\sevenCal=eusm7
\font\fiveCal=eusm5
\newfam\Calfam
  \textfont\Calfam=\tenCal
  \scriptfont\Calfam=\sevenCal
  \scriptscriptfont\Calfam=\fiveCal
\def\Cal{\fam\Calfam\tenCal}

\def\cC{{\Cal C}}
\def\cD{{\Cal D}}
\def\cE{{\Cal E}}
\def\cI{{\Cal I}}
\def\cK{{\Cal K}}

\def\cP{{\Cal P}}
\def\cX{{\Cal X}}

\def\det{{\mathop{\rm det}\nolimits}}

\def\Ricci{\mathop{\rm Ricci}\nolimits}

\def\dim{\mathop{\rm dim}\nolimits}
\def\codim{\mathop{\rm codim}\nolimits}

\def\Ker{\mathop{\rm Ker}\nolimits}
\def\Ricci{\mathop{\rm Ricci}\nolimits}

\def\dbar{\overline\partial}
\def\ddbar{\partial\overline\partial}
\def\\{\hfil\break}
\let\wt\widetilde

\let\ovl\overline

\def\hexnbr#1{\ifnum#1<10 \number#1\else
 \ifnum#1=10 A\else\ifnum#1=11 B\else\ifnum#1=12 C\else
 \ifnum#1=13 D\else\ifnum#1=14 E\else\ifnum#1=15 F\fi\fi\fi\fi\fi\fi\fi}
\def\msatype{\hexnbr\msafam}
\def\msbtype{\hexnbr\msbfam}
\mathchardef\compact="3\msatype62
\mathchardef\smallsetminus="2\msbtype72   \let\ssm\smallsetminus
\mathchardef\subsetneq="3\msbtype28
\def\buildo#1\over#2{\mathrel{\mathop{\null#2}\limits^{#1}}}
\def\buildu#1\under#2{\mathrel{\mathop{\null#2}\limits_{#1}}}

\def\square{{\hfill \hbox{
\vrule height 1.453ex  width 0.093ex  depth 0ex
\vrule height 1.5ex  width 1.3ex  depth -1.407ex\kern-0.1ex
\vrule height 1.453ex  width 0.093ex  depth 0ex\kern-1.35ex
\vrule height 0.093ex  width 1.3ex  depth 0ex}}}
\def\qed{\phantom{$\quad$}\hfill$\square$\medskip}

% Redefign item to align references on the left border
\let\Item\item
\def\item#1{\Item{$\rlap{\hbox{#1}}\kern\parindent\kern-5pt$}}

% Today's date
\def\today{\ifcase\month\or
January\or February\or March\or April\or May\or June\or July\or August\or
September\or October\or November\or December\fi \space\number\day,
\number\year}

%%%%%%%%%%%%%%%%%%%%%%%%%%%%%%%%%%%%%%%%%%%%%%%%%%%%%%%%%%%%%%%%%%%%%%%%%%%%%%

\null
\vskip 3cm
\centerline{\fourteenbf Numerical characterization of the}
\bigskip
\centerline{\fourteenbf K\"ahler cone of a compact K\"ahler manifold}

\null
\vskip 5pt
\centerline{\bf Jean-Pierre Demailly${}^*$ and Mihai Paun${}^{**}$}
\medskip
\noindent
$\strut\kern3.5cm$ \llap{${}^*$} Universit\'e Joseph Fourier, Grenoble\\
$\strut\kern3.5cm$ \llap{${}^{**}$} Universit\'e Louis Pasteur, Strasbourg

\null
\vskip 30pt
{\eightpoint

\noindent
{\bf Abstract.}
The goal of this work is give a precise numerical description of the
K\"ahler cone of a compact K\"ahler manifold. Our main result states
that the K\"ahler cone depends only on the intersection form of the
cohomology ring, the Hodge structure and the homology classes of
analytic cycles: if $X$ is a compact K\"ahler manifold, the K\"ahler
cone $\cK$ of $X$ is one of the connected components of the set
$\cP$ of real $(1,1)$ cohomology classes $\{\alpha\}$ which are
numerically positive on analytic cycles, i.e.\ $\int_Y\alpha^p>0$ for
every irreducible analytic set $Y$ in $X$, \hbox{$p=\dim Y$}.  This result is
new even in the case of projective manifolds, where it can be seen as
a generalization of the well-known Nakai-Moishezon criterion, and it
also extends previous results by Campana-Peternell and Eyssidieux.
The principal technical step is to show that every nef class $\{\alpha\}$
which has positive highest self-intersection number $\int_X\alpha^n>0$
contains a K\"ahler current; this is done by
using the Calabi-Yau theorem and a mass concentration technique for
Monge-Amp\`ere equations. The main result admits a number of variants and 
corollaries, including a description of the cone of numerically effective
$(1,1)$ classes and their dual cone. Another important consequence is 
the fact that for an arbitrary deformation $\cX\to S$ of compact 
K\"ahler manifolds, the K\"ahler cone of a very general fibre $X_t$ 
is ``independent'' of $t$, i.e.\ invariant by parallel transport under 
the $(1,1)$-component of the Gauss-Manin connection.
}
\vskip40pt

\noindent{\twelvebf 0. Introduction} 
\medskip

\noindent
The primary goal of this work is to study in great detail the structure of
the K\"ahler cone of a compact K\"ahler manifold. Our main result 
states that the K\"ahler cone depends only on the intersection form
of the cohomology ring, the Hodge structure and the homology classes of 
analytic cycles. More precisely, we have
\medskip

\noindent 
{\bf 0.1. Main Theorem.} {\sl Let $X$ be a compact K\"ahler manifold.
Then the K\"ahler cone $\cK$ of $X$ is one of the connected components
of the set $\cP$ of real $(1,1)$ cohomology classes $\{\alpha\}$
which are numerically positive on analytic cycles, i.e.\ $\int_Y\alpha^p>0$
for every irreducible analytic set $Y$ in $X$, $p=\dim Y$.
}
\medskip

This result is new even in the case of projective manifolds. It can be
seen as a generalization of the well-known Nakai-Moishezon criterion,
which provides a necessary and sufficient criterion for a line bundle
to be ample: {\sl a line bundle $L\rightarrow X$ on a projective
algebraic manifold $X$ is ample if and only if $\displaystyle
\int_Yc_1(L)^p> 0$, for every algebraic subset $Y\subset X$, $p=\dim
Y$.}  When $X$ is projective, we show that the condition
$\int_Y\alpha^p> 0$ characterizes precisely all K\"ahler classes, even
when $\{\alpha\}$ is not an integral class -- and even when
$\{\alpha\}$ lies outside the real Neron-Severi group
$NS_\bR(X)=NS(X)\otimes_\bZ\bR$. In the notation of the Main Theorem,
we have therefore $\cK=\cP$ if $X$ is projective.

This result extends a few special cases which had been proved earlier
by completely different methods: Campana-Peternell [CP90] showed 
that the Nakai-Moishezon criterion holds true for classes 
$\{\alpha\}\in NS_\bR(X)$. Quite recently, using $L^2$ cohomology 
techniques for infinite coverings of a projective algebraic manifold, 
P.~Eyssidieux [Eys00] obtained a version of the Nakai-Moishezon for 
all real combinations of $(1,1)$ cohomology classes which become integral 
after taking the pull-back to some finite or infinite covering.

The Main Theorem admits quite a number of useful variants and corollaries.
Two of them are descriptions of the cone of numerically 
effective (nef) classes -- see section~1 for the precise definition
of numerical effectivity on general compact complex manifolds.
\medskip

\noindent 
{\bf 0.2. Corollary.} {\sl Let $X$ be a compact
K\"ahler manifold. A $(1,1)$ cohomology class $\{\alpha\}$ on $X$
is nef if and only if there exists a K\"ahler metric $\omega$ on $X$
such that $\int_Y\alpha^k\wedge\omega^{p-k}\ge 0$ for all 
irreducible analytic sets $Y$ and all $k=1,2,\ldots,p=\dim Y$.
}
\medskip

\noindent 
{\bf 0.3. Corollary.} {\sl Let $X$ be a compact K\"ahler manifold. 
A $(1,1)$ cohomology class $\{\alpha\}$ on $X$ is nef if and only 
for every irreducible analytic set $Y$ in $X$, $p=\dim X$
and every K\"ahler metric $\omega$ on $X$ we have
$\int_Y\alpha\wedge\omega^{p-1}\ge 0$. In other words, the dual of the
nef cone is the closed convex cone generated by cohomology classes of 
currents of the form $[Y]\wedge\omega^{p-1}$ in $H^{n-1,n-1}(X,\bR)$.
}
\medskip

\noindent
We now briefly discuss the essential ideas involved in our approach. The
first basic result is a sufficient condition for a nef class to contain a
K\"ahler current. The proof is based on a technique of mass concentration
for Monge-Amp\`ere equations, using the Aubin-Calabi-Yau theorem [Yau78].
\medskip

\noindent 
{\bf 0.4. Theorem.} {\sl Let $(X,\omega)$ be a compact 
$n$-dimensional K\"ahler manifold and let $\{\alpha\}$ in $H^{1,1}(X,\bR)$
be a nef cohomology class such that $\int_X\alpha^n>0$. Then $\{\alpha\}$
contains a K\"ahler current $T$, that is, a closed positive current
$T$ such that $T\ge\delta\omega$ for some $\delta>0$. The current $T$
can be chosen smooth in the complement $X\ssm Z$ of an analytic
set, with logarithmic poles along~$Z$.
}
\medskip

In a first step, we show that the class $\{\alpha\}^p$ dominates a
small multiple of any $p$-codimensional analytic set $Y$ in $X$. As
we already mentioned, this is done by concentrating the mass on $Y$ in 
the Monge-Amp\`ere equation. We then apply this fact to the diagonal 
$\Delta\subset \wt X=X\times X$ to produce a closed positive current
$\Theta\in\{\pi_1^*\alpha+\pi_2^*\alpha\}^n$ which dominates 
$[\Delta]$ in $X\times X$. The desired K\"ahler current $T$ is 
easily obtained by taking a push-forward 
$\pi_{1*}(\Theta\wedge\pi_2^*\omega)$ of $\Theta$ to~$X$.

The technique produces a priori ``very singular'' currents, since we use a
weak compactness argument. However, we can apply the general
regularization theorem proved in [Dem92] to get a current which is
smooth outside an analytic set $Z$ and only has logarithmic poles
along $Z$.  The idea of using a Monge-Amp\`ere equation to force the
occurrence of positive Lelong numbers in the limit current was first
exploited in [Dem93], in the case when $Y$ is a finite set of points,
to get effective results for adjoints of ample line bundles (e.g.\ in
the direction of the Fujita conjecture). 

The use of higher dimensional subsets $Y$ in the mass concentration
process will be crucial here. However, the technical details are quite
different from the $0$-dimensional case used in [Dem93]; in fact, we
cannot rely any longer on the maximum principle, as in the case of
Monge-Amp\`ere equations with isolated Dirac masses in the right hand
side.  The new technique employed here is essentially taken from
[Pau00]; it was already proved there -- for projective manifolds --
that any big semi-positive $(1,1)$-class contains a K\"ahler
current. The Main Theorem is deduced from 0.4 by induction on
dimension, thanks to the following useful result which was already
observed in the second author's Thesis ([Pau98a, Pau98b]).  \medskip

\noindent 
{\bf 0.5. Proposition.} {\sl Let $X$ be a compact complex manifold $($or
complex space$)$. Then
\smallskip
\item{\rm (i)} The cohomology class of a closed positive 
$(1,1)$-current $\{T\}$ is nef if and only if the restriction 
$\{T\}_{|Z}$ is nef for every irreducible component $Z$ in the Lelong 
sublevel sets~$E_c(T)$. 
\smallskip
\item{\rm (ii)} The cohomology class of a K\"ahler current $\{T\}$ is 
a K\"ahler class $($i.e.\ the class of a {\rm smooth} K\"ahler form$)$ 
if and only if the restriction $\{T\}_{|Z}$ is a K\"ahler class for 
every irreducible component $Z$ in the Lelong sublevel sets~$E_c(T)$. 
\vskip0pt
}
\medskip

To derive the Main theorem from 0.4 and 0.5, it is enough to observe that 
any class $\{\alpha\}\in\ovl\cK\cap\cP$ is nef and such that
$\int_X\alpha>0$, therefore it contains a K\"ahler current. By the
induction hypothesis on dimension, $\{\alpha\}_{|Z}$ is K\"ahler
for all $Z\subset X$, hence $\{\alpha\}$ is a K\"ahler class on~$X$.

We want to stress that Theorem 0.4 is closely related to the
solution of the Grauert-Riemenschneider conjecture by Y.-T. Siu
([Siu85]); see also [Dem85] for a stronger result based on holomorphic
Morse inequalities, and T.~Bouche [Bou89], S.~Ji-B.~Shiffman [JS93],
L.~Bonavero [Bon93, Bon98] for other related results.  The results obtained
by Siu can be summarized as follows: {\sl Let $L$ be a hermitian
semi-positive line bundle on a compact $n$-dimensional complex
manifold~$X$, such that $\int_X c_1(L)^n>0$. Then $X$ is a Moishezon
manifold and $L$ is a big line bundle; the tensor powers of $L$ have a
lot of sections, $h^0(X, L^m)\geq Cm^n$ as $m\to+\infty$, and there
exists a singular hermitian metric on $L$ such that the curvature of
$L$ is positive, bounded away from $0$}. Again, Theorem 0.4 can be
seen as an extension of this result to non integral $(1,1)$
cohomology classes -- however, our proof only works so far for
K\"ahler manifolds, while the Grauert-Riemenschneider conjecture has
been proved on arbitrary compact complex manifolds.  In the same vein,
we prove the following result.
\medskip

\noindent 
{\bf 0.6. Theorem.} {\sl A compact complex manifold carries a K\"ahler
current if and only if it is bimeromorphic to a K\"ahler manifold
$($or equivalently, dominated by a K\"ahler manifold$)$.
}
\medskip

This class of manifolds is called the {\sl Fujiki class $\cC$}. If we 
compare this result with the solution of the Grauert-Riemenschneider
conjecture, it is tempting to make the following conjecture which would
somehow encompass both results.
\medskip

\noindent 
{\bf 0.7. Conjecture.} {\sl Let $X$ be a compact complex manifold of
dimension~$n$. Assume that $X$ possesses a nef cohomology class $\{\alpha\}$
of type $(1,1)$ such that $\int_X\alpha^n>0$. Then $X$ is in the
Fujiki class~$\cC$.\\
$($Also, $\{\alpha\}$ would contain a K\"ahler current, as it
follows from Theorem 0.4 if Conjecture 0.7 is proved$)$.
}
\medskip

We want mention here that most of the above results were already known in 
the cases of complex surfaces (i.e.\ dimension $2$),
thanks to the work of N.~Buchdahl [Buc99, 00] and A.~Lamari [Lam99a, 99b] --
it turns out that there exists a very neat characterization of nef 
classes on arbitrary surfaces -- K\"ahler or not.

The main Theorem has an important application to the deformation theory of
compact K\"ahler manifolds, which we prove in Section 5.
\medskip

\noindent
{\bf 0.8. Theorem.} {\sl Let $\cX\to S$ be a deformation of compact K\"ahler
manifolds over an irreducible base~$S$. Then there exists a countable union 
$S'=\bigcup S_\nu$ of analytic subsets $S_\nu\subsetneq S$, such that 
the K\"ahler cones
$\cK_t\subset H^{1,1}(X_t,\bC)$ are invariant over $S\ssm S'$ under 
parallel transport with respect to the $(1,1)$-projection $\nabla^{1,1}$ of
the Gauss-Manin connection.}
\medskip

We moreover conjecture (see 5.2 for details) that the K\"ahler property is
open with respect to the countable Zariski topology on the base $S$ of
a deformation of arbitrary compact complex manifolds.

Shortly after this work was completed, Daniel Huybrechts [Huy01] informed
us that our Main Theorem can be used to calculate the K\"ahler cone of
a very general hyperk\"ahler manifold: the K\"ahler cone is then equal
to one of the connected components of the positive cone defined by the
Beauville-Bogomolov quadratic form.  This closes the gap in his
original proof of the projectivity criterion for hyperk\"ahler manifolds 
([Huy99], Theorem 3.11).

We are grateful to Arnaud Beauville, Christophe Mourougane and
Philippe Eyssidieux for helpful discussions, which were part of the 
motivation for looking at the questions investigated here.

\vfill\eject

\noindent {\twelvebf 1. Nef cohomology classes and K\"ahler currents}
\medskip
Let $X$ be a complex analytic manifold. Throughout this paper, we
denote by $n$ the complex dimension $\dim_\bC X$. As is well known, a K\"ahler
metric on $X$ is a smooth {\sl real} form of type $(1,1)$
$$
\omega(z)=i\sum_{1\le j,k\le n}\omega_{jk}(z)dz_j\wedge d\ovl z_k,
$$
that is, $\ovl\omega=\omega$ or equivalently
$\ovl{\omega_{jk}(z)}=\omega_{kj}(z)$, such that
\smallskip
{\parindent=1.12cm
\noindent
\item{$(1.1')$\phantom{'}} 
$\omega(z)$ is positive definite at every point [$(\omega_{jk}(z))$
is a positive definite hermitian matrix];
\smallskip\noindent
\item{$(1.1'')$} $d\omega=0$ when $\omega$ is viewed as a real $2$-form,
i.e.\ $\omega$ is symplectic.
\smallskip
}
One says that $X$ is K\"ahler (or is of K\"ahler type) if $X$ possesses
a K\"ahler metric $\omega$.
To every closed real (resp.\ complex) valued $k$-form $\alpha$ we 
associate its De Rham cohomology class $\{\alpha\}\in H^k(X,\bR)$
(resp.\ $\{\alpha\}\in H^k(X,\bC)$), and to every $\dbar$-closed
form $\alpha$ of pure type $(p,q)$ we associate its Dolbeault
cohomology class $\{\alpha\}\in H^{p,q}(X,\bC)$. On a compact K\"ahler
manifold we have a canonical Hodge decomposition
$$
H^k(X,\bC)=\bigoplus_{p+q=k}H^{p,q}(X,\bC).\leqno(1.2)
$$
In this work, we are especially interested in studying the {\sl K\"ahler cone}
$$
\cK\subset H^{1,1}(X,\bR):=H^{1,1}(X,\bC)\cap H^2(X,\bR),\leqno(1.3)
$$
which is by definition the set of cohomology classes $\{\omega\}$ of all
$(1,1)$-forms associated with K\"ahler metrics. Clearly, $\cK$ is an
open convex cone in $H^{1,1}(X,\bR)$, since a small perturbation of
a K\"ahler form is still a K\"ahler form. The closure $\ovl\cK$ of
the K\"ahler cone is equally important. Since we want to possibly
consider non K\"ahler manifolds, we have to consider ``$\ddbar$-cohomology'' 
groups
$$
H^{p,q}_{\ddbar}(X,\bC):=\{\hbox{$d$-closed $(p,q)$-forms}\}/
\ddbar\{\hbox{$(p-1,q-1)$-forms}\}.\leqno(1.4)
$$
When $(X,\omega)$ is compact K\"ahler, it is well known (from the so-called
$\ddbar$-lemma) that we have an isomorphism 
$H_{\ddbar}^{p,q}(X,\bC)\simeq H^{p,q}(X,\bC)$ with the more usual
Dolbeault groups. Notice that there are always canonical morphisms
$$
H_{\ddbar}^{p,q}(X,\bC)\to H^{p,q}(X,\bC), \qquad
H_{\ddbar}^{p,q}(X,\bC)\to H^{p+q}_{\rm DR}(X,\bC)
$$
($\ddbar$-cohomology is ``more precise'' than Dolbeault or De Rham cohomology).
This allows us to define numerically effective classes in a fairly
general situation (see also [Dem90b, Dem92, DPS94]).
\bigskip

\noindent {\sl {\bf 1.5. Definition.} Let $X$ be a compact complex manifold 
equipped with a hermitian positive $($non necessarily K\"ahler$)$ 
metric~$\omega$. A class $\{\alpha\}\in H_{\ddbar}^{1,1}(X,\bR)$
is said to be numerically effective
$($or nef for brevity$)$ if for every $\varepsilon>0$ there is a
representative $\alpha_\varepsilon=\alpha+i\ddbar\varphi_\varepsilon
\in\{\alpha\}$ such that $\alpha_\varepsilon\ge-\varepsilon\omega$.
}
\bigskip

\noindent
If $(X,\omega)$ is compact K\"ahler, a class $\{\alpha\}$ is nef if and 
only if $\{\alpha+\varepsilon\omega\}$ is a K\"ahler class for every 
$\varepsilon>0$,
i.e., a class $\{\alpha\}\in H^{1,1}(X,\bR)$ is nef if and only if it belongs 
to the closure $\ovl\cK$ of the K\"ahler cone. (Also, if $X$ is 
projective algebraic, a divisor $D$ is nef in the sense of algebraic 
geometers, that is, $D\cdot C\ge 0$ for every irreducible curve $C\subset X$,
if and only if $\{D\}\in\ovl\cK$, so the definitions fit together; see
[Dem90b, Dem92] for more details).

In the sequel, we will make a heavy use of currents, especially the theory 
of closed positive currents. Recall that a current $T$ is a differential
form with distribution coefficients. In the complex situation, we are
interested in currents
$$
T=i^{pq}\sum_{|I|=p,|J|=q} T_{I,J}\,dz_I\wedge d\ovl z_J\qquad
\hbox{($T_{I,J}$ distributions on $X$),}
$$
of pure bidegree $(p,q)$, with $dz_I=dz_{i_1}\wedge\ldots\wedge dz_{i_p}$ 
as usual. We say that $T$ is positive if $p=q$ and
$\sum \lambda_I\ovl\lambda_JT_{I,J}$ is a positive distribution (i.e.\
a positive measure) for all possible choices of complex 
coefficients~$\lambda_I$, $|I|=p$. Alternatively, the space of 
$(p,q)$-currents can be seen as the dual space of the Fr\'echet space of
smooth $(n-p,n-q)$-forms, and $(n-p,n-q)$ is called the bidimension of~$T$.
By Lelong [Lel57], to every analytic set $Y\subset X$ of codimension $p$
is associated a current $T=[Y]$ defined by
$$
\langle[Y], u\rangle=\int_Y u,\qquad u\in\cD_{n-p,n-p}(X),
$$
and $[Y]$ is a closed positive current of bidegree $(p,p)$ and
bidimension \hbox{$(n-p,n-p)$}. The theory of positive currents can be
easily extended to complex spaces $X$ with singularities; one then
simply defines the space of currents to be the dual of space of smooth
forms, defined as forms on the regular part $X_{\rm reg}$ which,
near $X_{\rm sing}$, locally extend as smooth forms on an open set of
$\bC^N$ in which $X$ is locally embedded (see e.g.\ [Dem85] for more 
details).
\bigskip

\noindent {\sl {\bf 1.6. Definition.} A K\"ahler current on a compact
complex space $X$ is a closed positive current $T$ of bidegree $(1,1)$
which satisfies $T\ge\varepsilon\omega$ for some $\varepsilon>0$ and
some smooth positive hermitian form $\omega$ on $X$.
} \medskip

\noindent
When $X$ is a (non singular) compact complex manifold, we consider 
the {\sl pseudo-effective} cone $\cE\subset H^{1,1}_{\ddbar}(X,\bR)$, defined
as the set of $\ddbar$-cohomology classes of closed positive 
$(1,1)$-currents. By the weak compactness of bounded sets in the 
space of currents, this is always a closed (convex) cone. When $X$ is 
K\"ahler, a K\"ahler current is just an element of the interior $\cE^\circ$
of $\cE$, and we have
$$
\cK\subset\cE^\circ.
$$
The inclusion may be strict, however, even when $X$ is K\"ahler, and the 
existence of a K\"ahler current on $X$ does not necessarily imply that 
$X$ admits a (smooth) K\"ahler form, as we will see in section 3 (and 
therefore $X$ need not be a K\"ahler manifold$\,$!). 
\vskip 30pt 

\noindent {\twelvebf 2. Concentration of mass for big nef classes}
\medskip

In this section, we show in full generality that every big and nef
cohomology class on a compact K\"ahler manifold contains a K\"ahler
current. The proof is based on a mass concentration technique for
Monge-Amp\`ere equations, using the Aubin-Calabi-Yau theorem.  We
first start by an easy lemma, which was (more or less) already
observed in [Dem90a]. Recall that a {\sl quasi-plurisubharmonic}
function $\psi$, by definition, is a function which is locally the sum
of a plurisubharmonic function and of a smooth function, or
equivalently, a function such that $i\ddbar\psi$ is locally bounded
below by a negative smooth $(1,1)$-form.
\medskip

\noindent 
{\bf 2.1. Lemma.} {\sl Let $X$ be a compact complex manifold~$X$ 
equipped with a K\"ahler metric 
$\omega=i\sum_{1\le j,k\le n}\omega_{jk}(z)dz_j\wedge d\ovl z_k$
and let $Y\subset X$ be an analytic subset of~$X$. Then there exist
globally defined quasi-plurisubharmonic potentials $\psi$ and 
$(\psi_\varepsilon)_{\varepsilon\in{}]0,1]}$ on $X$,
satisfying the following properties.
\smallskip
\item {\rm (i)} The function $\psi$ is smooth on $X\ssm Y$, satisfies
$i\ddbar\psi\ge -A\omega$ for some $A>0$, and $\psi$
has logarithmic poles along~$Y$, i.e., locally near $Y$
$$
\psi(z)\sim \log\sum_k|g_k(z)|+O(1)
$$
where $(g_k)$ is a local system of generators of the ideal sheaf $\cI_Y$
of $Y$ in $X$.
\smallskip
\item {\rm(ii)} We have 
$\psi=\lim_{\varepsilon\to 0}\downarrow\psi_\varepsilon$ 
and the $\psi_\varepsilon$ possess a uniform Hessian estimate
$$
i\ddbar\psi_\varepsilon\ge -A\omega\qquad\hbox{on~$X$}.
$$
\smallskip
\item {\rm(iii)} Consider the family of hermitian metrics
$$
\omega_\varepsilon := \omega + {1\over 2A} i\ddbar\psi_\varepsilon
\ge {1\over 2}\omega.
$$
For any point $x_0\in Y$ and any neighborhood $U$ of $x_0$, 
the volume element of $\omega_\varepsilon$ has
a uniform lower bound 
$$
\int_{U\cap V_\varepsilon}\omega_\varepsilon^n\ge\delta(U)>0,
$$
where $V_\varepsilon=\{z\in X\,;\;\psi(z)<\log\varepsilon\}$
is the ``tubular neighborhood'' of radius $\varepsilon$ around $Y$.
\smallskip
\item {\rm(iv)} For every integer $p\ge 0$, the family of positive 
currents $\omega_\varepsilon^p$ is bounded in mass. Moreover,
if $Y$ contains an irreducible component $Y'$ of codimension $p$,
there is a uniform lower bound
$$
\int_{U\cap V_\varepsilon}\omega_\varepsilon^p\wedge\omega^{n-p}\ge
\delta_p(U)>0
$$
in any neighborhood $U$ of a regular point $x_0\in Y'$. In particular, 
any weak limit $\Theta$ of $\omega_\varepsilon^p$ as $\varepsilon$
tends to $0$ satisfies $\Theta\ge\delta'[Y']$ for some
$\delta'>0$.
\vskip0pt } 
\medskip

\noindent
{\sl Proof.} By compactness of $X$, there is a covering of $X$ by open
coordinate balls~$B_j$, $1\le j\le N$, such that $\cI_Y$ is generated by 
finitely many holomorphic functions $(g_{j,k})_{1\le k\le m_j}$ on a 
neighborhood of $\ovl B_j$. We take a partition of unity
$(\theta_j)$ subordinate to $(B_j)$ such that $\sum\theta_j^2=1$ on $X$,
and define
$$
\eqalign{
\psi(z)&={1\over 2}\log\sum_j\theta_j(z)^2\sum_k|g_{j,k}(z)|^2,\cr
\psi_\varepsilon(z)&={1\over 2}\log(e^{2\psi(z)}+\varepsilon^2)
={1\over 2}\log\Big(\sum_{j,k}\theta_j(z)^2|g_{j,k}(z)|^2+\varepsilon^2
\Big).\cr}
$$
Moreover, we consider the family of $(1,0)$-forms with support in $B_j$
such that
$$
\gamma_{j,k}=\theta_j\partial g_{j,k}+2g_{j,k}\partial\theta_j.
$$
Straightforward calculations yield
$$
\leqalignno{
\dbar\psi_\varepsilon&={1\over 2}
{\sum_{j,k}\theta_jg_{j,k}\ovl{\gamma_{j,k}}\over e^{2\psi}+
\varepsilon^2},\cr
i\ddbar\psi_\varepsilon&={i\over 2}\left(
{\sum_{j,k}\gamma_{j,k}\wedge\ovl{\gamma_{j,k}}\over e^{2\psi}+
\varepsilon^2}-
{\sum_{j,k}\theta_j\ovl{g_{j,k}}\gamma_{j,k}
\wedge \sum_{j,k}\theta_jg_{j,k}\ovl{\gamma_{j,k}}
\over (e^{2\psi}+\varepsilon^2)^2}\right),&(2.2)\cr
&\qquad{}+i\;
{\sum_{j,k}|g_{j,k}|^2(\theta_j\ddbar\theta_j-\partial\theta_j\wedge
\dbar\theta_j)\over e^{2\psi}+\varepsilon^2}.\cr
}
$$
As $e^{2\psi}=\sum_{j,k}\theta_j^2|g_{j,k}|^2$, the first big sum in
$i\ddbar\psi_\varepsilon$ is nonnegative by the Cauchy-Schwarz inequality;
when viewed as a hermitian form, the value of this sum on a tangent 
vector $\xi\in T_X$ is simply
$$
{1\over 2}\left(
{\sum_{j,k}|\gamma_{j,k}(\xi)|^2\over e^{2\psi}+\varepsilon^2}-
{\big|\sum_{j,k}\theta_j\ovl{g_{j,k}}\gamma_{j,k}(\xi)\big|^2\over
(e^{2\psi}+\varepsilon^2)^2}\right)\ge
{1\over 2}{\varepsilon^2\over(e^{2\psi}+\varepsilon^2)^2}
\sum_{j,k}|\gamma_{j,k}(\xi)|^2.\leqno(2.3)
$$
Now, the second sum involving $\theta_j\ddbar\theta_j-\partial\theta_j\wedge
\dbar\theta_j$ in (2.2) is uniformly bounded below by a fixed negative 
hermitian form $-A\omega$, $A\gg 0$, and therefore
\hbox{$i\ddbar\psi_\varepsilon\ge -A\omega$}. Actually, for every pair
of indices $(j,j')$ we have a bound
$$
C^{-1}\le \sum_k|g_{j,k}(z)|^2/\sum_k|g_{j',k}(z)|^2\le C\qquad
\hbox{on $\ovl B_j\cap\ovl B_{j'}$},
$$
since the generators $(g_{j,k})$ can be expressed as holomorphic
linear combinations of the $(g_{j',k})$ by Cartan's theorem A (and
vice versa). It follows easily that all terms $|g_{j,k}|^2$ are uniformly
bounded by $e^{2\psi}+\varepsilon^2$. In particular, $\psi$ and
$\psi_\varepsilon$ are quasi-plurisubharmonic, and we see that 
(i) and (ii) hold true. By construction, 
the real $(1,1)$-form $\omega_\varepsilon:=\omega+{1\over 2A}
i\ddbar\psi_\varepsilon$ satisfies $\omega_\varepsilon\ge {1\over 2}\omega$,
hence it is K\"ahler and its eigenvalues with respect to $\omega$ are at 
least equal to~$1/2$. 

Assume now that we are in a neighborhood $U$ of a regular point
$x_0\in Y$ where $Y$ has codimension~$p$.  Then
\hbox{$\gamma_{j,k}=\theta_j\partial g_{j,k}$} at $x_0$, hence the
rank of the system of $(1,0)$-forms $(\gamma_{j,k})_{k\ge 1}$ is at
least equal to $p$ in a neighborhood of~$x_0$. Fix a holomorphic
locate coordinate system $(z_1,\ldots,z_n)$ such that 
\hbox{$Y=\{z_1=\ldots=z_p=0\}$} near $x_0$, and let $S\subset T_X$ 
be the holomorphic subbundle generated by $\partial/\partial z_1,\ldots,
\partial/\partial z_p$. This choice ensures that the rank of the system
of $(1,0)$-forms $(\gamma_{j,k|S})$ is everywhere equal to~$p$.
By (1,3) and the minimax principle applied to the $p$-dimensional subspace
$S_z\subset T_{X,z}$, we see that the $p$-largest eigenvalues of 
$\omega_\varepsilon$ are bounded below by 
$c\varepsilon^2/(e^{2\psi}+\varepsilon^2)^2$.
\smallskip
However, we can even restrict the form defined in (2.3) to the 
$(p-1)$-dimen\-sional subspace $S\cap\Ker\tau$ where 
$\tau(\xi):=\sum_{j,k}\theta_j\ovl{g_{j,k}}\gamma_{j,k}(\xi)$,
to see that the $(p-1)$-largest eigenvalues of $\omega_\varepsilon$
are bounded below by $c/(e^{2\psi}+\varepsilon^2)$, $c>0$. The
$p$-th eigenvalue is then bounded by 
$c\varepsilon^2/(e^{2\psi}+\varepsilon^2)^2$ and the
remaining $(n-p)$-ones by~$1/2$. From this we infer
$$
\eqalign{
\omega_\varepsilon^n
&\ge c\,{\varepsilon^2\over(e^{2\psi}+\varepsilon^2)^{p+1}}\omega^n
\qquad\hbox{near $x_0$},\cr
\omega_\varepsilon^p
&\ge c\,{\varepsilon^2\over(e^{2\psi}+\varepsilon^2)^{p+1}}\Big
(i\sum_{1\le\ell\le p}\gamma_{j,k_\ell}\wedge
\ovl{\gamma_{j,k_\ell}}\Big)^p
\cr}
$$
where $(\gamma_{j,k_\ell})_{1\le \ell\le p}$ is a suitable $p$-tuple
extracted from the $(\gamma_{j,k})$, such that\break 
$\bigcap_\ell\Ker\gamma_{j,k_\ell}$ is a smooth complex
(but not necessarily holomorphic) subbundle of codimension $p$ of $T_X$;
by the definition of the forms $\gamma_{j,k}$, this
subbundle must coincide with $T_Y$ along $Y$. From this, properties
(iii) and (iv) follow easily; actually, up to constants,
we have $e^{2\psi}+\varepsilon^2\sim|z_1|^2+\ldots+|z_p|^2+\varepsilon^2$
and
$$
i\sum_{1\le\ell\le p}\gamma_{j,k_\ell}\wedge
\ovl{\gamma_{j,k_\ell}}\ge c\,i\ddbar(|z_1|^2+\ldots+|z_p|^2)-
O(\varepsilon)i\ddbar|z|^2\qquad\hbox{on $U\cap V_\varepsilon$},
$$
hence, by a straightforward calculation,
$$
\omega_\varepsilon^p\wedge\omega^{n-p}\ge
c\big(i\ddbar\log(|z_1|^2+\ldots+|z_p|^2+\varepsilon^2)\big)^p\wedge
\big(i\ddbar(|z_{p+1}|^2+\ldots+|z_n|^2)\big)^{n-p}
$$
on $U\cap V_\varepsilon$; notice also that
$\omega_\varepsilon^n\ge 2^{-(n-p)}\omega_\varepsilon^p\wedge\omega^{n-p}$,
so any lower bound for the volume of 
$\omega_\varepsilon^p\wedge\omega^{n-p}$ will also produce a bound
for the volume of $\omega_\varepsilon^n$. As it is well known, the $(p,p)$-form
$$
\Big({i\over 2\pi}\ddbar\log(|z_1|^2+\ldots+|z_p|^2+\varepsilon^2)\Big)^p
\qquad\hbox{on $\bC^n$}
$$
can be viewed as the pull-back to $\bC^n=\bC^p\times\bC^{n-p}$ of 
the Fubini-Study volume form of the complex $p$-dimensional
projective space of dimension $p$ containing $\bC^p$ as an affine
Zariski open set, rescaled by the dilation ratio~$\varepsilon$. Hence
it converges weakly to the current of integration on the $p$-codimensional
subspace $z_1=\ldots=z_p=0$. Moreover the volume contained in any
compact tubular cylinder
$$
\{|z'|\le C\varepsilon\}\times K''\subset \bC^p\times\bC^{n-p}
$$
depends only on $C$ and $K$ (as one sees after rescaling by $\varepsilon$).
The fact that $\omega_\varepsilon^p$ is uniformly bounded in mass can
be seen easily from the fact that
$$
\int_X \omega_\varepsilon^p\wedge\omega^{n-p}=\int_X\omega^n,
$$
as $\omega$ and $\omega_\varepsilon$ are in the same K\"ahler class.
Let $\Theta$ be any weak limit of $\omega_\varepsilon^p$. By what we
have just seen, $\Theta$ carries non zero mass on every 
$p$-codimensional component $Y'$ of $Y$, for instance near every
regular point. However, standard results of the theory of currents
(support theorem and Skoda's extension result) imply that
${\bf 1}_{Y'}\Theta$ is a closed positive current and that
${\bf 1}_{Y'}\Theta=\lambda[Y']$ is a nonnegative multiple
of the current of integration on $Y'$. The fact that the mass
of $\Theta$ on $Y'$ is positive yields~$\lambda>0$. 
Lemma~2.1 is proved.\qed
\bigskip

\noindent {\bf 2.4. Remark.} In fact, we did not really make use of the
fact that $\omega$ is K\"ahler. Lemma~2.1 would still be true without
this assumption. The only difficulty would be to show that 
$\omega_\varepsilon^p$ is still locally bounded in mass when $\omega$
is an arbitrary hermitian metric. This can be done by using a resolution 
of singularities which converts $\cI_Y$ into an invertible sheaf defined 
by a divisor with normal crossings -- and by doing some standard
explicit calculations. As we do not need the more general form of
Lemma~2.1, we will omit these technicalities. 
\bigskip

\noindent Let us now recall the following very deep result concerning 
Monge-Amp\`ere equations on compact K\"ahler manifolds (see~[Yau78]).
\medskip

\noindent 
{\bf 2.5. Theorem {\rm (Yau.)}} {\sl Let $(X,\omega)$ 
be a compact K\"ahler manifold and $n=\dim X$. Then for any smooth
volume form $f>0$ such that $\int_X f = \int_X\omega^n$,
there exist a K\"ahler metric $\widetilde \omega = \omega+i\ddbar\varphi$
in the same K\"ahler class as $\omega$, such that $\widetilde \omega^n=f.$}
\medskip

\noindent
In other words, one can prescribe the volume form $f$ of the K\"ahler metric
$\wt\omega\in\{\omega\}$, provided that the total volume $\int_Xf$ is
equal to the expected value $\int_X\omega^n$. Since the Ricci curvature form
of $\wt\omega$ is $\Ricci(\wt\omega):=
-{i\over2\pi}\ddbar\log\det(\wt\omega)=-{i\over2\pi}\ddbar\log f$,
this is the same as prescribing the curvature form $\Ricci(\wt\omega)=\rho$,
given any $(1,1)$-form $\rho$ representing $c_1(X)$. Using this, we prove
\medskip

\noindent 
{\bf 2.6. Proposition.} {\sl Let $(X,\omega)$ be a compact 
$n$-dimensional K\"ahler manifold and let $\{\alpha\}$ in $H^{1,1}(X,\bR)$
be a nef cohomology class
such that $\alpha^n>0$ $($a ``big'' class$)$. For every
$p$-codimensional analytic set $Y\subset X$, there exists a closed positive
current $\Theta\in\{\alpha\}^p$ of bidegree $(p,p)$ such that
$\Theta\ge\delta[Y]$ for some $\delta>0$.
}
\medskip

\noindent
{\sl Proof.} Let us associate with $Y$ a family $\omega_\varepsilon$
of K\"ahler metrics as in Lemma~2.1. The 
class $\{\alpha+\varepsilon\omega\}$ is a K\"ahler
class, so by Yau's theorem we can find a representative
$\alpha_\varepsilon=\alpha+\varepsilon\omega+i\ddbar\varphi_\varepsilon$
such that
$$
\alpha_\varepsilon^n=C_\varepsilon\omega_\varepsilon^n,\leqno(2.7)
$$
where
$$
C_\varepsilon = {\int_X\alpha_\varepsilon^n\over\int_X\omega_\varepsilon^n}
={\int_X(\alpha+\varepsilon\omega)^n\over\int_X\omega^n}\ge
C_0={\int_X\alpha^n\over\int_X\omega^n}>0.
$$
Let us denote by 
$$
\lambda_1(z)\le\ldots\le\lambda_n(z)
$$
the eigenvalues of $\alpha_\varepsilon(z)$ with respect to
$\omega_\varepsilon(z)$, at every point $z\in X$ (these functions
are continuous with respect to $z$, and of course depend also 
on~$\varepsilon$). The equation (2.7) is equivalent to the fact that
$$
\lambda_1(z)\ldots\lambda_n(z)=C_\varepsilon\leqno(2.7')
$$
is constant, and the most important observation for us is that the
constant $C_\varepsilon$ is bounded away from $0$, thanks to
our assumption $\int_X\alpha^n>0$.

Fix a regular point $x_0\in Y$ and a small neighborhood $U$ (meeting only
the irreducible component of $x_0$ in $Y$). By Lemma 2.1, we have a
uniform lower bound
$$
\int_{U\cap V_\varepsilon}\omega_\varepsilon^p\wedge\omega^{n-p}
\ge\delta_p(U)>0.\leqno(2.8)
$$
Now, by looking at the $p$ smallest (resp.\ $(n-p)$ largest) eigenvalues
$\lambda_j$ of $\alpha_\varepsilon$ with respect to $\omega_\varepsilon$, 
we find
$$
\leqalignno{
\alpha_\varepsilon^p&\ge \lambda_1\ldots\lambda_p\,\omega_\varepsilon^p,
&(2.9')\cr
\alpha_\varepsilon^{n-p}\wedge\omega_\varepsilon^p
&\ge {1\over n!}\lambda_{p+1}\ldots\lambda_n\,
\omega_\varepsilon^n,&(2.9'')\cr
}
$$
The last inequality $(2.9'')$ implies
$$
\int_X\lambda_{p+1}\ldots\lambda_n\,\omega_\varepsilon^n\le
n!\int_X\alpha_\varepsilon^{n-p}\wedge\omega_\varepsilon^p=
n!\int_X(\alpha+\varepsilon\omega)^{n-p}\wedge\omega^p\le M
$$
for some constant $M>0$ (we assume $\varepsilon\le 1$, say). 
In particular, for every $\delta>0$, the subset 
$E_\delta\subset X$ of points $z$ such that 
$\lambda_{p+1}(z)\ldots\lambda_n(z)>M/\delta$ satisfies
$\int_{E_\delta}\omega_\varepsilon^n\le\delta$, hence
$$
\int_{E_\delta}\omega_\varepsilon^p\wedge\omega^{n-p}\le 2^{n-p}
\int_{E_\delta}\omega_\varepsilon^n\le 2^{n-p}\delta.\leqno(2.10)
$$
The combination of (2.8) and (2.10) yields
$$
\int_{(U\cap V_\varepsilon)\ssm E_\delta}\omega_\varepsilon^p\wedge\omega^{n-p}
\ge \delta_p(U)-2^{n-p}\delta.
$$
On the other hand $(2.7')$ and $(2.9')$ imply
$$
\alpha_\varepsilon^p\ge{C_\varepsilon\over\lambda_{p+1}\ldots\lambda_n}
\omega_\varepsilon^p
\ge{C_\varepsilon\over M/\delta}\omega_\varepsilon^p\qquad
\hbox{on $(U\cap V_\varepsilon)\ssm E_\delta$}.
$$
From this we infer
$$
\int_{U\cap V_\varepsilon}\alpha_\varepsilon^p\wedge\omega^{n-p}\ge
{C_\varepsilon\over M/\delta}
\int_{(U\cap V_\varepsilon)\ssm E_\delta}
\omega_\varepsilon^p\wedge\omega^{n-p}\ge
{C_\varepsilon\over M/\delta}(\delta_p(U)-2^{n-p}\delta)>0
\leqno(2.11)
$$
provided that $\delta$ is taken small enough, e.g.\ 
$\delta=2^{-(n-p+1)}\delta_p(U)$. The family of $(p,p)$-forms
$\alpha_\varepsilon^p$ is uniformly bounded in mass since
$$
\int_X\alpha_\varepsilon^p\wedge\omega^{n-p}
=\int_X(\alpha+\varepsilon\omega)^p\wedge\omega^{n-p}\le\hbox{Const}.
$$
Inequality $(2.11)$ implies that any weak limit $\Theta$ of 
$(\alpha_\varepsilon^p)$
carries a positive mass on~$U\cap Y$. By Skoda's extension theorem [Sk81], 
${\bf 1}_Y\Theta$ is a closed positive current with support in $Y$, 
hence ${\bf 1}_Y\Theta=\sum\lambda_j[Y_j]$ is a combination of the 
various components $Y_j$ of $Y$ with coefficients $\lambda_j>0$. 
Our construction shows that $\Theta$ belongs to the cohomology class 
$\{\alpha\}^p$. Proposition 2.6 is proved.\qed
\bigskip

\noindent
We can now prove the main result of this section.
\medskip

\noindent 
{\bf 2.12. Theorem.} {\sl Let $(X,\omega)$ be a compact 
$n$-dimensional K\"ahler manifold and let $\{\alpha\}$ in $H^{1,1}(X,\bR)$
be a nef cohomology class such that $\int_X\alpha^n>0$. Then $\{\alpha\}$
contains a K\"ahler current $T$, that is, a closed positive current
$T$ such that $T\ge\delta\omega$ for some $\delta>0$.
}
\medskip

\noindent
{\sl Proof.} The trick is to apply Proposition 2.6 to the diagonal
$\wt Y=\Delta$ in the product manifold $\wt X=X\times X$.
Let us denote by $\pi_1$ and $\pi_2$ the two projections of
$\wt X=X\times X$ onto $X$. It is clear that $\wt X$ admits
$$
\wt\omega=\pi_1^*\omega+ \pi_2^*\omega
$$
as a K\"ahler metric, and that the class of
$$
\wt\alpha=\pi_1^*\alpha + \pi_2^*\alpha
$$
is a nef class on $\wt X$ [it is a limit of the K\"ahler classes
$\pi_1^*(\alpha+\varepsilon\omega) + \pi_2^*(\alpha+\varepsilon\omega)$].
Moreover, by Newton's binomial formula
$$
\int_{X\times X}\wt\alpha^{2n}={2n \choose n}\Big(\int_X\alpha^n\Big)^2>0.
$$
The diagonal is of codimension $n$ in $\wt X$, hence by Proposition 2.6
there exists a closed positive $(n,n)$-current 
$\Theta\in\{\wt\alpha^n\}$ such that
$\Theta\ge\varepsilon[\Delta]$ for some $\varepsilon>0$. We define
the $(1,1)$-current $T$ to be the push-forward
$$
T=c\,\pi_{1*}(\Theta\wedge\pi_2^*\omega)
$$
for a suitable constant $c>0$ which will be determined later. By the
lower estimate on~$\Theta$, we have
$$
T\ge c\varepsilon\,\pi_{1*}([\Delta]\wedge\pi_2^*\omega)=
c\varepsilon\,\omega,
$$
thus $T$ is a K\"ahler current. On the other hand, as 
$\Theta\in\{\wt\alpha^n\}$, the current $T$ belongs to the
cohomology class of the $(1,1)$-form
$$
c\,\pi_{1*}(\wt\alpha^n\wedge\pi_2^*\omega)(x)=c\int_{y\in Y}
\big(\alpha(x)+\alpha(y)\big)^n\wedge\omega(y),
$$
obtained by a partial integration in $y$ with respect to $(x,y)\in X\times X$.
By Newton's binomial formula again, we see that
$$
c\,\pi_{1*}(\wt\alpha^n\wedge\pi_2^*\omega)(x)=c\Big(\int_Xn\alpha(y)^{n-1}\wedge
\omega(y)\Big)\alpha(x)
$$
is proportional to $\alpha$. Therefore, we need only take 
$c=\big(\int_Xn\alpha^{n-1}\wedge
\omega\big)^{-1}$ to ensure that $T\in\{\alpha\}$. Notice that
$\alpha$ is nef and $\{\alpha\}\le C\{\omega\}$ for sufficiently large $C>0$,
so we have
$$
\int_X\alpha^{n-1}\wedge\omega\ge{1\over C}\int_X\alpha^n>0.
$$
Theorem 2.12 is proved.\qed
\vskip 30pt

\noindent{\twelvebf 3. Regularization theorems for K\"ahler currents} 
\medskip

It is not true that a K\"ahler current can be regularized to produce
a smooth K\"ahler metric. However, by the general regularization 
theorem for closed currents proved in [Dem92] (cf.\ Proposition 3.7),
it can be regularized up to some logarithmic poles along analytic subsets.

Before stating the result, we need a few preliminaries. If $T$ is a 
closed positive current on a compact complex manifold $X$, we can write
$$
T = \alpha+i\ddbar\psi\leqno(3.1)
$$
where $\alpha$ is a global smooth closed $(1,1)$-form on $X$, and $\psi$
a quasi-plurisubharmonic function on $X$. To see this (cf.\ also [Dem92])
take an open covering of $X$ by open coordinate balls $B_j$ and 
plurisubharmonic potentials 
$\psi_j$ such that $T=i\ddbar\psi_j$ on $B_j$. Then, if $(\theta_j)$ is
a partition of unity subordinate to $(B_j)$, it is easy to see that
$\psi=\sum \theta_j\psi_j$ is quasi-plurisubharmonic and that
$\alpha:=T-i\ddbar\psi$ is smooth (so that
$i\ddbar\psi=T-\alpha\ge -\alpha$). For any other decomposition
$T=\alpha'+i\ddbar\psi'$ as in (3.1), we have $\alpha'-\alpha=
-i\ddbar(\psi'-\psi)$, hence $\psi'-\psi$ is smooth.
\bigskip

\noindent {\bf 3.2. Regularization theorem.} {\sl
Let $X$ be a compact complex manifold equipped with a hermitian metric
$\omega$. Let $T=\alpha+i\ddbar\psi$ be a closed $(1,1)$-current on~$X$,
where $\alpha$ is smooth and $\psi$ is a quasi-plurisubharmonic function. 
Assume that $T\ge\gamma$ for some real $(1,1)$-form $\gamma$ on $X$ with 
real coefficients. Then there exists 
a sequence $T_k=\alpha+i\ddbar\psi_k$ of closed $(1,1)$-currents such that
\smallskip
\item{\rm (i)} $\psi_k$ $($and thus $T_k)$ is smooth on the complement
$X\ssm Z_k$ of an analytic set $Z_k$, and the $Z_k$'s form an increasing
sequence
$$
Z_0\subset Z_1\subset\ldots\subset Z_k\subset\ldots \subset X.
$$
\smallskip
\item{\rm (ii)} There is a uniform estimate $T_k\ge \gamma-\delta_k\omega$
with $\lim\downarrow\delta_k=0$ as $k$ tends to $+\infty$.
\smallskip
\item{\rm (iii)} The sequence $(\psi_k)$ is non increasing, and we 
have $\lim\downarrow\psi_k=\psi$. As a consequence, $T_k$ converges
weakly to $T$ as  $k$ tends to $+\infty$.
\smallskip
\item{\rm (iv)} Near $Z_k$, the potential $\psi_k$ has logarithmic poles,
namely, for every $x_0\in Z_k$, there is a neighborhood $U$ of $x_0$
such that $\psi_k(z)=\lambda_k\log\sum_\ell|g_{k,\ell}|^2+O(1)$ for 
suitable holomorphic functions $(g_{k,\ell})$ on $U$ and $\lambda_k>0$.
Moreover, there is a $($global$)$ proper modification $\mu_k:\wt X_k\to X$ 
of $X$, obtained as a sequence of blow-ups with smooth centers, such that 
$\psi_k\circ\mu_k$ can be written locally on $\wt X_k$ as
$$
\psi_k\circ\mu_k(w)=\lambda_k\big(\sum n_\ell\log|\wt g_\ell|^2+f(w)\big)
$$
where $(\wt g_\ell=0)$ are local generators of suitable $($global$)$
divisors $D_\ell$ on $\wt X_k$ such that $\sum D_\ell$ has normal crossings,
$n_\ell$ are positive integers, and the $f$'s are smooth functions 
on~$\wt X_k$.
\vskip0pt
}
\medskip

\noindent
{\sl Sketch of the proof.} We briefly indicate the main ideas, since the
proof can only be reconstructed by patching together arguments which
appeared in different places (although the core the proof is entirely
in [Dem92]). After replacing $T$ with $T-\alpha$, we can assume that 
$\alpha=0$ and
$T=i\ddbar\psi\ge\gamma$. Given a small $\varepsilon>0$, we select a 
covering of $X$ by open balls $B_j$ together with holomorphic coordinates 
$(z^{(j)})$ and real numbers $\beta_j$ such that
$$
0\le \gamma-\beta_j\,i\ddbar|z^{(j)}|^2\le \varepsilon\,i\ddbar|z^{(j)}|^2
\qquad\hbox{on $B_j$}
$$
(this can be achieved just by continuity of $\gamma$, after 
diagonalizing $\gamma$ at the center of the balls). We now take a partition
of unity $(\theta_j)$ subordinate to $(B_j)$ such that $\sum\theta_j^2=1$,
and define
$$
\psi_k(z)={1\over 2k}\log\sum_j\theta_j^2
e^{2k\beta_j|z^{(j)}|^2}\sum_{\ell\in\bN}|g_{j,k,\ell}|^2
$$
where $(g_{j,k,\ell})$ is a Hilbert basis of the Hilbert space of holomorphic
functions $f$ on $B_j$ such that
$$
\int_{B_j}|f|^2e^{-2k(\psi-\beta_j|z^{(j)}|^2)}<+\infty.
$$
Notice that by the Hessian estimate $i\ddbar\psi\ge\gamma\ge\beta_j\,
i\ddbar|z^{(j)}|^2$, the weight involved in the $L^2$ norm is
plurisubharmonic. It then follows from the proof of Proposition 3.7
in [Dem92] that all properties (i)--(iv) hold true, except possibly the
fact that the sequence $\psi_k$ can be chosen to be non increasing, and
the existence of the modification in (iv). However, the multiplier ideal 
sheaves of the weights \hbox{$k(\psi-\beta_j|z^{(j)}|^2)$} are generated by
the $(g_{j,k,\ell})_\ell$ on $B_j$, and these sheaves glue together into
a global coherent multiplier ideal sheaf $\cI(k\psi)$ on $X$
(see [DEL99]); the modification $\mu_k$ is then obtained by blowing-up 
the ideal sheaf $\cI(k\psi)$ so that $\mu_k^*\cI(k\psi)$ is an invertible
ideal sheaf associated with a normal crossing divisor (Hironaka [Hir63]).
The fact that $\psi_k$ can be chosen to be non increasing follows from 
a quantitative version of the ``subadditivity of multiplier ideal sheaves'' 
which is proved in Step 3 of the proof of Theorem 2.2.1 in [DPS00]
(see also ([DEL99]). (Anyway, this property will not be used here, so 
the reader may wish to skip the details).\qed
\bigskip

For later purposes, we state the following useful results, which are
borrowed essentially from the PhD thesis of the second author.
\bigskip

\noindent {\bf 3.3. Proposition {\rm [Pau98a, Pau98b]}.} {\sl Let $X$ be 
a compact complex space and let $\{\alpha\}$ be a $\ddbar$-cohomology
class of type $(1,1)$ on~$X$ $($where $\alpha$ is a smooth 
represen\-ta\-tive$)$.
\smallskip
\item{\rm (i)} If the restriction $\{\alpha\}_{|Y}$ to an analytic subset
$Y\subset X$ is K\"ahler on $Y$, there exists a smooth representative
$\alpha'=\alpha+i\ddbar\varphi$ which is K\"ahler on a neighborhood $U$
of~$Y$.
\smallskip
\item{\rm (ii)} If the restrictions $\{\alpha\}_{|Y_1}$, $\{\alpha\}_{|Y_2}$
to any pair of analytic subsets $Y_1,\,Y_2\subset X$ are nef 
$($resp.\ K\"ahler$)$,
then $\{\alpha\}_{|Y_1\cup Y_2}$ is nef $($resp.\ K\"ahler$)$.
\smallskip
\item{\rm (iii)} Assume that $\{\alpha\}$ contains a K\"ahler current $T$
and that the restriction $\{\alpha\}_{|Y}$ to every irreducible component 
$Y$ in the Lelong sublevel sets $E_c(T)$ is a K\"ahler class on $Y$.  
Then $\{\alpha\}$ is a K\"ahler class on~$X$.
\smallskip
\item{\rm (iv)} Assume that $\{\alpha\}$ contains a closed positive
$(1,1)$-current $T$ and that the restriction $\{\alpha\}_{|Y}$ to every 
irreducible component $Y$ in the Lelong sublevel sets
$E_c(T)$ is nef on $Y$. Then $\{\alpha\}$ is nef on~$X$.
\vskip0pt
}
\medskip

By definition, $E_c(T)$ is the set of points $z\in X$ such that the
Lelong number $\nu(T,z)$ is at least equal to $c$ (for given $c>0$).
A deep theorem of Siu ([Siu74]) asserts that all $E_c(T)$ are analytic
subsets of $X$. Notice that the concept of a $\ddbar$-cohomology class
is well defined on an arbitrary complex space (although many of the
standard results on De Rham or Dolbeault cohomology of non singular
spaces will fail for singular spaces!). The concepts of K\"ahler
classes and nef classes are still well defined [a~K\"ahler form on a
singular space $X$ is a $(1,1)$-form which is locally bounded below by
the restriction of a smooth positive $(1,1)$-form in a non singular
ambient space for $X$, and a nef class is a class contained
representatives bounded below by $-\varepsilon\omega$ for every
$\varepsilon>0$, where $\omega$ is a smooth positive $(1,1)$-form].
\medskip

\noindent
{\sl Sketch of the proof.} (i) We can assume that $\alpha_{|Y}$ itself
is a K\"ahler form. If $Y$ is smooth, we simply take $\psi$ to
be equal to a large constant times the square of the hermitian distance 
to~$Y$. This will produce positive eigenvalues in $\alpha+i\ddbar\psi$
along the normal directions of $Y$, while the eigenvalues are already
postive on $Y$. When $Y$ is singular, we just use the same argument
with respect to a stratification of $Y$ by smooth manifolds, and an
induction on the dimension of the strata ($\psi$ can be left untouched
on the lower dimensional strata).
\smallskip

\noindent
(ii) Let us first treat the K\"ahler case. By (i), there are smooth functions 
$\varphi_1$, $\varphi_2$ on $X$ such that $\alpha+i\ddbar\varphi_j$ is 
K\"ahler on a neighborhood $U_j$ of $Y_j$, $j=1,2$. Also, by Lemma 2.1,
there exists a quasi plurisubharmonic function $\psi$ on $X$
which has logarithmic poles on $Y_1\cap Y_2$ and is smooth on $X\ssm(Y_1\cap Y_2)$. We define
$$
\varphi=\wt{\max}(\varphi_1+\delta\psi,\varphi_2-C)
$$
where $\delta\ll 1$, $C\gg 1$ are constants and $\wt{\max}$ is a 
regularized max function. Then $\alpha+i\ddbar\varphi$ is K\"ahler on 
$U_1\cap U_2$. Moreover, for $C$ large, $\varphi$ coincides with 
$\varphi_1+\delta\psi$ on $Y_1\ssm U_2$ and with
$\varphi_2-C$ on a small neighborhood of $W$ of $Y_1\cap Y_2$. 
Take smaller neighborhoods $U_1'\compact U_1$, $U'_2\compact U_2$ such
that $U'_1\cap U'_2\subset W$. We can extend $\varphi_{|U'_1\cap U_2}$ 
to a neighborhood $V$ of $Y_1\cup Y_2$ by taking 
$\varphi=\varphi_1+\delta\psi$ on a neighborhood of
$Y_1\ssm U_2$ and $\varphi=\varphi_2-C$ on $U'_2$. The use of a cut-off
function equal to $1$ on a neighborhood of $V'\compact V$ of $Y_1\cup Y_2$
finally allows us to get a function $\varphi$ defined everywhere on $X$,
such that $\alpha+i\ddbar\varphi$ is K\"ahler on a neighborhood of 
$Y_1\cup Y_2$ (if $\delta$ is small enough). The nef case is entirely 
similar, except that we deal with currents $T$ such that 
$T\ge-\varepsilon\omega$ instead of K\"ahler currents.
\smallskip

\noindent
(iii) By the regularization theorem 3.2, we may assume that the
singularities of the K\"ahler current $T=\alpha+i\ddbar\psi$ are just
logarithmic poles (since $T\ge\gamma$ with $\gamma$ positive definite,
the small loss of positivity resulting from 3.2~(ii) still yields a
K\"ahler current~$T_k$). Hence $\psi$ is smooth on $X\ssm Z$ for a suitable
analytic set $Z$ which, by construction, is contained in $E_c(T)$ for
$c>0$ small enough. We use (i), (ii) and the hypothesis that $\{T\}_{|Y}$
is K\"ahler for every component $Y$ of $Z$ to get a neighborhood $U$
of $Z$ and a smooth potential $\varphi_U$ on $X$ such that
$\alpha+i\ddbar\varphi_U$ is K\"ahler on $U$. Then the smooth potential
equal to the regularized maximum $\varphi=\wt{\max}(\psi,\varphi_U-C)$
produces a K\"ahler form $\alpha+i\ddbar\varphi$ on $X$ for
$C$ large enough (since we can achieve $\varphi=\psi$ on $X\ssm U$).
The nef case (iv) is entirely similar.\qed
\bigskip

\noindent
{\bf 3.4. Theorem. } {\sl A compact complex manifold $X$ admits a
K\"ahler current if and only if it is bimeromorphic to a K\"ahler manifold,
or equivalently, if it admits a proper K\"ahler modification.
$($The class of such manifolds is the so-called {\rm Fujiki class $\cC$)}.}
\medskip

\noindent
{\sl Proof.} If $X$ is bimeromorphic to a K\"ahler manifold $Y$, 
Hironaka's desingularization theorem implies that there exists
a blow-up $\widetilde Y$ of $Y$ (obtained by a sequence of blow-ups with
smooth centers) such that the bimeromorphic map from
$Y$ to $X$ can be resolved into a modification $\mu:\widetilde Y\to X$. Then
$\widetilde Y$ is K\"ahler and the push-forward $T=\mu_*\wt\omega$
of a K\"ahler form $\wt\omega$ on $\widetilde Y$ provides a K\"ahler 
current on~$X$. In fact, if $\omega$ is a smooth hermitian form on $X$,
there is a constant $C$ such that $\mu^*\omega\le C\wt\omega$ (by compactness
of $\wt Y$), hence
$$
T=\mu_*\wt\omega\ge \mu_*(C^{-1}\mu^*\omega)=C^{-1}\omega.
$$
Conversely, assume that $X$ admits a K\"ahler current $T$. By Theorem
3.2~(iv), there exists a K\"ahler current $T'=T_k$ ($k\gg1$) in the same 
$\ddbar$-cohomology class as $T$, and a modification
$\mu:\wt X\to X$ such that
$$
\mu^*T'=\lambda[\wt D]+\wt\alpha\qquad\hbox{on $\wt X$},
$$
where $\wt D$ is a divisor with normal crossings, $\wt\alpha$
a smooth closed $(1,1)$-form and $\lambda>0$. The form $\wt\alpha$
must be semi-positive, more
precisely we have $\wt\alpha\ge\varepsilon\mu^*\omega$ as soon as
$T'\ge\varepsilon\omega$. This is not enough to produce a K\"ahler
form on $\wt X$ (but we are not very far...). Suppose that $\wt X$ 
is obtained as a tower of blow-ups
$$
\wt X=X_N\to X_{N-1}\to\cdots \to X_1\to X_0=X,
$$
where $X_{j+1}$ is the blow-up of $X_j$ along a smooth center $Y_j\subset X_j$.
Denote by $E_{j+1}\subset X_{j+1}$ the exceptional divisor, and let
$\mu_j:X_{j+1}\to X_j$ be the blow-up map. Now, we use the following
simple
\smallskip

\noindent{\bf 3.5. Lemma.} {\sl For every K\"ahler current 
$T_j$ on $X_j$, there exists $\varepsilon_{j+1}>0$ and a smooth form $u_{j+1}$
in the $\ddbar$-cohomology class of $[E_{j+1}]$ such that
$$
T_{j+1}=\mu_j^\star T_j-\varepsilon_{j+1}u_{j+1}
$$ 
is a K\"ahler current on~$X_{j+1}$.}
\smallskip

\noindent The reason is that the line bundle
${\cal O}(-E_{j+1})|E_{j+1}$ is equal to ${\cal O}_{P(N_j)}(1)$ where $N_j$
is the normal bundle to $Y_j$ in $X_j$, hence it is positive along the
normal directions to $Y_j\,$; assume furthermore that $T_j\ge\delta_j\omega_j$
for suitable $0<\delta_j\ll 1$ and a hermitian form $\omega_j$ on $X_j\,$;
then
$$
\mu_j^\star T_j-\varepsilon_{j+1}u_{j+1}\ge
\delta_j\mu_j^\star\omega_j-\varepsilon_{j+1}u_{j+1}
$$
where $\mu_j^*\omega_j$ is semi-positive on $X_{j+1}$, positive definite on 
$X_{j+1}\ssm E_{j+1}$, and also positive definite along the ``horizontal''
directions of $Y_j$ on $E_{j+1}$. The statement is then easily proved
by taking $\varepsilon_{j+1}\ll\delta_j$ and by using a compactness
argument.

If $\widetilde u_j$ is the pull-back of $u_j$ to the final blow-up
$\widetilde X$, we conclude inductively that $\mu^\star
T'-\sum\varepsilon_j\widetilde u_j$ is a K\"ahler current. Therefore
the smooth form
$$
\wt\omega:=\wt\alpha-\sum\varepsilon_j\widetilde u_j
=\mu^\star T'-\sum\varepsilon_j\widetilde u_j-\lambda[D]
$$
is K\"ahler and we see that $\wt X$ is a K\"ahler manifold.\qed
\bigskip

\noindent{\bf 3.6. Remark.} A special case of Theorem 3.4 is the following
characterization of Moishezon varieties (i.e.\ manifolds which are
bimeromorphic to projective algebraic varieties or, equivalently, whose
algebraic dimension is equal to their complex dimension): 
\smallskip
\noindent
{\sl A compact complex manifold $X$ is Moishezon if and only if $X$ 
possesses a K\"ahler current $T$ such that the De Rham cohomology class
$\{T\}$ is rational, i.e.\ $\{T\}\in H^2(X,\bQ)$}. 
\smallskip
\noindent
In fact, in the above proof, we get an integral current $T$ if we take
the push forward $T=\mu_*\wt\omega$ of an integral ample class 
$\{\wt\omega\}$ on $Y$, where $\mu:Y\to X$ is a projective model 
of~$Y$. Conversely, if $\{T\}$ is rational, we can take
the $\varepsilon_j's$ to be rational in Lemma~3.5. This produces
at the end a K\"ahler metric $\wt\omega$ with rational De Rham cohomology
class on $\wt X$. Therefore $\wt X$ is projective by the Kodaira
embedding theorem. This result was already observed in [JS93]
(see also [Bon93, Bon98] for a more general perspective based
on a singular version of holomorphic Morse inequalities).
\vskip30pt

\noindent{\twelvebf 4. Numerical characterization of the K\"ahler cone} 
\medskip

We are now in a good position to derive what we consider to be the 
main result of this work.
\bigskip

\noindent
{\bf 4.1. Main theorem.} {\sl Let $X$ be a compact K\"ahler manifold, and let
$$
\cP\subset H^{1,1}(X,\bR)
$$ 
be the set of $(1,1)$-classes $\{\alpha\}$ which are ``numerically positive''
an analytic cycles:
$$
\int_Y\alpha^p>0
$$
for all irreducible analytic subsets $Y\subset X$, $\dim Y=p$. Then the
K\"ahler cone $\cK$ is one of the connected components of $\cP$ $($and
also, one of the connected components of the interior $\cP^\circ$ 
of~$\cP)$.
}
\medskip

\noindent
{\sl Proof.} By definition $\cK$ is open, and clearly $\cK\subset\cP$
(thus $\cK\subset\cP^\circ$). We claim that $\cK$ is also closed in~$\cP$.
In fact, consider a class $\{\alpha\}\in\ovl\cK\cap\cP$. This means
that $\{\alpha\}$ is a nef class which satisfies all numerical conditions
defining~$\cP$. Let $Y\subset X$ be an arbitrary analytic subset. We prove 
by induction on $\dim Y$ that $\{\alpha\}_{|Y}$ is K\"ahler. If $Y$ has 
several components, Proposition 3.3~(ii) reduces the situation to the
case of the irreducible components of~$Y$, so we may assume that $Y$ is 
irreducible. Let $\mu:\wt Y\to Y$ be a desingularization of $Y$, obtained
via a finite sequence of blow-ups with smooth centers in $X$. Then
$\wt Y$ is a smooth K\"ahler manifold and $\{\mu^*\alpha\}$ is a nef class
such that
$$
\int_{\wt Y}(\mu^*\alpha)^p=\int_Y\alpha^p>0,\qquad p=\dim Y.
$$
By Theorem 2.12, there exists a K\"ahler current $\wt T$ on $\wt Y$ which
belong to the class $\{\mu^\star\alpha\}$. Then $T:=\mu_*\wt T$ is
a K\"ahler current on $Y$, contained in the class $\{\alpha\}$.
By the induction hypothesis, the class $\{\alpha\}_{|Z}$ is K\"ahler
for every irreducible component $Z$ of $E_c(T)$ (since $\dim Z\le p-1$).
Proposition 3.3~(iii) now shows that $\{\alpha\}$ is K\"ahler on~$Y$.
In the case $Y=X$, we get that $\{\alpha\}$ itself is K\"ahler, hence
$\{\alpha\}\in\cK$ and $\cK$ is closed in~$\cP$. This implies that
$\cK$ is a union of connected components of~$\cP$. However, since $\cK$ is
convex, it is certainly connected, and only one component can be
contained in~$\cK$.\qed
\bigskip

\noindent
{\bf 4.2. Remark.} In all examples that we are aware of, the cone
$\cP$ is open. Moreover, Theorem 4.1 shows that the connected
component of any K\"ahler class $\{\omega\}$ in $\cP$ is open in
$H^{1,1}(X,\bR)$ (actually an open convex cone...).  However, it might
still happen that $\cP$ carries some boundary points on the other
components. It turns out that there exist examples for which $\cP$ is
not connected. Let us consider for instance a complex torus 
$X=\bC^n/\Lambda$. It is well-known that a generic torus $X$ does not
possess any analytic subset except finite subsets and $X$ itself. 
In that case, the numerical positivity is expressed by the single
condition $\int_X\alpha^n>0$. However, on a torus, $(1,1)$-classes are 
in one-to-one correspondence with constant hermitian forms $\alpha$ on
$\bC^n$. Thus, for $X$ generic, $\cP$ is the set of hermitian forms on 
$\bC^n$ such that $\det(\alpha)>0$, and Theorem 4.1 just expresses the 
elementary result of linear algebra saying that the set $\cK$ of 
positive definite forms is one of the connected components of the open set
$\cP=\{\det(\alpha)>0\}$ of hermitian forms of positive determinant
(the other components, of course, are the sets of forms of signature 
$(p,q)$, $p+q=n$, $q$ even).
\bigskip

One of the drawbacks of Theorem 4.1 is that the characterization of the
K\"ahler cone still involves the choice of an undetermined connected component.
However, it is trivial to derive the following (weaker) variants, which do
not involve the choice of a connected component.
\bigskip

\noindent
{\bf 4.3. Theorem.} {\sl Let $(X,\omega)$ be a compact K\"ahler manifold
and let $\{\alpha\}$ be a $(1,1)$ cohomology class in~$H^{1,1}(X,\bR)$.
The following properties are equivalent.
\smallskip
\item{\rm (i)} $\{\alpha\}$ is K\"ahler.
\smallskip
\item{\rm (ii)} For every irreducible analytic set $Y\subset X$, $\dim Y=p$, 
and every $t\ge 0$
$$
\int_Y(\alpha+t\omega)^p>0.
$$
\smallskip
\item{\rm (iii)} For every irreducible analytic set $Y\subset X$, $\dim Y=p$,
$$
\int_Y\alpha^k\wedge\omega^{p-k}>0\qquad\hbox{for $k=1,\ldots,p$}.
$$
}
\medskip

\noindent
{\sl Proof.} It is obvious that (i)$\,\Rightarrow\,$(iii)$\,\Rightarrow\,$(ii),
so we only need to show that (ii)$\,\Rightarrow\,$(i). Assume that
condition (ii) holds true. For $t_0$ large enough, $\alpha+t_0\omega$ 
is a K\"ahler class.
The segment $(\alpha+t_0\omega)_{t\in[0,t_0]}$ is a connected set intersecting
$\cK$ which is contained in $\cP$, thus it is entirely 
contained in $\cK$ by Theorem 4.1. We infer that $\{\alpha\}\in\cK$, as
desired.\qed
\bigskip

\noindent
{\bf 4.4. Theorem.} {\sl Let $X$ be a compact K\"ahler manifold and
let $\{\alpha\}\in H^{1,1}(X,\bR)$ be a  $(1,1)$ cohomology class.
The following properties are equivalent.
\smallskip
\item{\rm(i)} $\{\alpha\}$ is nef.
\smallskip
\item{\rm(ii)} There exists a K\"ahler class $\omega$ such that
$$
\int_Y\alpha^k\wedge\omega^{p-k}\ge 0
$$
for every irreducible analytic set $Y\subset X$, $\dim Y=p$, and every
$k=1,2,\ldots, p$.
\smallskip
\item{\rm(iii)} For every irreducible analytic set $Y\subset X$, $\dim Y=p$, 
and every K\"ahler class $\{\omega\}$ on $X$
$$
\int_Y\alpha\wedge\omega^{p-1}\ge 0.
$$
}
\medskip

\noindent
{\sl Proof.} Clearly (i)$\,\Rightarrow\,$(ii) and (i)$\,\Rightarrow\,$(iii).
\smallskip
\noindent
(ii)$\,\Rightarrow\,$(i). If $\{\alpha\}$ satisfies the inequalities in 
(ii), then the class $\{\alpha+\varepsilon\omega\}$ satisfies the 
corresponding  strict inequalities for every $\varepsilon>0$. Therefore
$\{\alpha+\varepsilon\omega\}$ is K\"ahler by Theorem 4.3, and 
$\{\alpha\}$ is nef.
\smallskip

\noindent
(iii)$\,\Rightarrow\,$(i). This is the most tricky part. For every 
integer $p\ge 1$, there exists a polynomial identity of the form
$$
(y-\delta x)^p - (1-\delta)^px^p
= (y-x)\int_0^1 A_p(t,\delta)\big((1-t)x+ty\big)^{p-1}\,dt
\leqno(4.5)
$$
where $A_p(t,\delta)=\sum_{0\le m\le p}a_m(t)\delta^m\in
\bQ[t,\delta]$ is a polynomial of degree${}\le p-1$ in~$t$ (moreover, 
the polynomial $A_p$ is unique under this limitation for the  degree). 
To see this, we observe that $(y-\delta x)^p - (1-\delta)^px^p$
vanishes identically for $x=y$, so it is divisible by $y-x$. By
homogeneity in $(x,y)$, we have an expansion of the form
$$
(y-\delta x)^p - (1-\delta)^px^p = (y-x)
\sum_{0\le\ell\le p-1,\,0\le m\le p} b_{\ell,m}x^\ell y^{p-1-\ell}\delta^m
$$
in the ring $\bZ[x,y,\delta]$. Formula (4.5) is then equivalent to
$$
b_{\ell,m} = \int_0^1 a_m(t){p-1\choose\ell}(1-t)^\ell t^{p-1-\ell}\,dt.
\leqno(4.5')
$$
Since $(U,V)\mapsto\int_0^1U(t)V(t)dt$ is a non degenerate linear pairing
on the space of polynomials of degree${}\le p-1$ and since
$({p-1\choose\ell}(1-t)^\ell t^{p-1-\ell})_{0\le\ell\le p-1}$ is a basis
of this space, $(4.5')$ can be achieved for a unique choice of the
polynomials $a_m(t)$. A~straightforward calculation shows that $A_p(t,0)=1$
identically. We can therefore choose $\delta_0\in[0,1[$ so small that 
$A_p(t,\delta)>0$ for all $t\in[0,1]$, $\delta\in[0,\delta_0]$ and 
$p=1,2,\ldots,n$.

Now, fix a K\"ahler metric $\omega$ such that $\omega'=\alpha+\omega$ 
is K\"ahler (if necessary, multiply $\omega$ by a large constant to reach 
this). A substitution $x=\omega$ and $y=\omega'$ in our polynomial
identity yields
$$
(\alpha+(1-\delta)\omega)^p-(1-\delta)^p\omega^p=\int_0^1
A_p(t,\delta)\,\alpha\wedge\big((1-t)\omega+t\omega'\big)^{p-1}dt.
$$
For every irreducible analytic subset $Y\subset X$ of dimension $p$ we find
$$
\eqalign{
\int_Y(\alpha+(1-\delta)\omega)^p&{}-(1-\delta)^p\int_Y\omega^p\cr
&=\int_0^1A_p(t,\delta)dt\Big(
\int_Y\alpha\wedge\big((1-t)\omega+t\omega'\big)^{p-1}\Big).\cr
}
$$
However, $(1-t)\omega+t\omega'$ is K\"ahler and therefore
$\int_Y\alpha\wedge\big((1-t)\omega+t\omega'\big)^{p-1}\ge 0$ by 
condition (iii). This implies $\int_Y(\alpha+(1-\delta)\omega)^p>0$
for all $\delta\in[0,\delta_0]$. We have produced a segment entirely
contained in $\cP$ such that one extremity $\{\alpha+\omega\}$ is in~$\cK$,
so the other extremity $\{\alpha+(1-\delta_0)\omega\}$ is also in $\cK$.
By repeating the argument inductively, we see that 
$\{\alpha+(1-\delta_0)^\nu\omega\}\in\cK$ for every integer $\nu\ge 0$.
From this we infer that $\{\alpha\}$ is nef, as desired.\qed
\bigskip

Since condition 4.4~(iii) is linear with respect to $\alpha$, we can also
view this fact as a characterization of the dual cone of the nef cone,
in the space of real cohomology classes of type $(n-1,n-1)$. We can state
\medskip

\noindent
{\bf 4.6. Proposition.} {\sl Let $X$ be a compact $n$-dimensional manifold.
The dual cone of the nef cone $\ovl\cK\subset H^{1,1}(X,\bR)$ in 
$H^{n-1,n-1}(X,\bR)$ is the closed convex cone generated by cohomology
classes of bidimension $(1,1)$ currents of the form 
\hbox{$[Y]\wedge\omega^{p-1}$}, $p=\dim Y$, where $Y$ runs over the 
collection of irreducible analytic subsets of $X$ and $\{\omega\}$ over
the set of K\"ahler classes of~$X$.
}
\bigskip

In the case of projective manifolds, we get stronger and simpler
versions of the above statements. All these can be seen as an extension 
of the Nakai-Moishezon criterion to arbitrary $(1,1)$-classes classes 
(not just integral $(1,1)$-classes as in the usual Nakai-Moishezon
criterion). Apart from the special cases already mentioned in the
introduction ([CP90], [Eys00]), these results seem to be entirely new.
\medskip

\noindent
{\bf 4.7. Theorem.} {\sl Let $X$ be a projective algebraic manifold.
Then $\cK=\cP$. Moreover, we have the following numerical characterizations.
\smallskip
\item{\rm(i)} a $(1,1)$-class $\{\alpha\}\in H^{1,1}(X,\bR)$ 
is K\"ahler if and only if $\int_Y\alpha^p>0$ 
for every irreducible analytic set $Y\subset X$, $p=\dim Y$.
\smallskip
\item{\rm(ii)} a $(1,1)$-class $\{\alpha\}\in H^{1,1}(X,\bR)$ 
is nef if and only if $\int_Y\alpha^p\ge 0$ 
for every irreducible analytic set $Y\subset X$, $p=\dim Y$.
\smallskip
\item{\rm(iii)} a $(1,1)$-class $\{\alpha\}\in H^{1,1}(X,\bR)$ 
is nef if and only if $\int_Y\alpha\wedge\omega^{p-1}\ge 0$ 
for every irreducible analytic set $Y\subset X$, $p=\dim Y$, and
every K\"ahler class $\{\omega\}$ on $X$.
\vskip0pt
}
\medskip

\noindent
{\sl Proof.} (i) We take $\omega=c_1(A,h)$ equal to the curvature form of a 
very ample line bundle $A$ on $X$, and we apply the numerical conditions 
as they are expressed in 4.3~(ii). For every $p$-dimensional algebraic 
subset $Y$ in $X$ we have
$$
\int_Y\alpha^k\wedge\omega^{p-k}=\int_{Y\cap H_1\cap\ldots \cap H_{p-k}}
\omega^k
$$
for a suitable generic complete intersection $Y\cap H_1\cap\ldots \cap H_{p-k}$
of $Y$ by members of the linear system~$|A|$. This shows that $\cP=\cK$.
\smallskip
\noindent
(ii) The nef case follows by considering $\alpha+\varepsilon\omega$, and
letting $\varepsilon>0$ tend to~$0$.
\smallskip
\noindent
(iii) is true more generally for any compact K\"ahler manifold.\qed
\bigskip

\noindent
{\bf 4.8. Remark.} In the case of a divisor $D$ (i.e., of an integral
class $\{\alpha\}$) on a projective algebraic manifold $X$, it is well
know that $\{\alpha\}$ is nef if and only if $D\cdot C=\int_C\alpha\ge 0$ for
every algebraic curve $C$ in $X$. This result completely fails when
$\{\alpha\}$ is not an integral class -- this is the same as saying that
the dual cone of the nef cone, in general, is bigger than the closed
convex cone generated by cohomology classes of effective curves.
Any surface such that the Picard number $\rho$ is less than $h^{1,1}$
provides a counterexample (any generic abelian surface or any
generic projective K3 surface is thus a counterexample). In particular,
in 4.7~(iii), it is not sufficient to merely consider the {\sl integral}
K\"ahler classes $\{\omega\}$.
\vskip30pt

\noindent{\twelvebf 5. Deformations of compact K\"ahler manifolds} 
\medskip

Let $\pi:\cX\to S$ be a deformation of non singular compact K\"ahler manifolds,
i.e.\ a proper analytic map between reduced complex spaces, with smooth
K\"ahler fibres, such that
the map is is trivial fibration locally near every point of $\cX$ (this is
of course the case if $\pi:\cX\to S$ is smooth, but here we do not want 
to require $S$ to be smooth; however we will always assume $S$ to be 
irreducible -- hence connected as well). 

We wish to investigate the behaviour of the K\"ahler cones $\cK_t$ of
the various fibres $X_t=\pi^{-1}(t)$, as $t$ runs over $S$. Because of the
assumption of local triviality of $\pi$, the topology of $X_t$ is
locally constant, and therefore so are the cohomology groups
$H^k(X_t,\bC)$. Each of these forms a locally constant vector bundle 
over $S$, whose associated sheaf of sections is the direct image sheaf
$R^k\pi_*(\bC_\cX)$. This locally constant system of $\bC$-vector space
contains as a sublattice the locally constant system of
integral lattices $R^k\pi_*(\bZ_\cX)$. As a consequence, the Hodge
bundle $t\mapsto H^k(X_t,\bC)$ carries a natural flat connection
$\nabla$ which is known as the {\it Gauss-Manin connection}. 

Thanks to D.~Barlet's theory of cycle spaces [Bar75], one can attach to every 
reduced complex space $X$ a reduced cycle space $C^p(X)$ parametrizing 
its compact analytic cycles of a given complex dimension $p$. In our 
situation, there is a relative cycle space $C^p(\cX/S)\subset C^p(\cX)$
which consists of all cycles contained in the fibres of $\pi:X\to S$.
It is equipped with a canonical holomorphic projection 
$$
\pi_p:C^p(\cX/S)\to S.
$$
Moreover, as the fibres $X_t$ are K\"ahler, it is known that the
restriction of $\pi_p$ to the connected components of $C^p(\cX/S)$
are proper maps. Also, there is a cohomology class (or degree) map 
$$
C^p(\cX/S) \to R^{2q}\pi_*(\bZ_\cX),\qquad Z\mapsto\{[Z]\}
$$
commuting with the projection to $S$, which to every compact
analytic cycle $Z$ in $X_t$ associates its cohomology class 
$\{[Z]\}\in H^{2q}(X_t,\bZ)$,
where $q=\codim Z=\break\dim X_t-p$. Again by the K\"ahler property 
(bounds on volume and Bishop compactness theorem), the map
$C^p(\cX/S) \to R^{2q}\pi_*(\bZ_\cX)$ is proper.

As is well known, the Hodge filtration
$$
F^p(H^k(X_t,\bC))=\bigoplus_{r+s=k,r\ge p}H^{r,s}(X_t,\bC)
$$
defines a {\sl holomorphic} subbundle of $H^k(X_t,\bC)$ (with respect to its
locally constant structure). On the other hand, the Dolbeault groups are 
given by
$$
H^{p,q}(X_t,\bC) = F^p(H^k(X_t,\bC))\cap\overline{F^{k-p}(H^k(X_t,\bC))},\qquad
k=p+q,
$$
and they form {\sl real analytic} subbundles of $H^k(X_t,\bC)$. We are 
interested especially in the decomposition
$$
H^2(X_t,\bC)=H^{2,0}(X_t,\bC)\oplus H^{1,1}(X_t,\bC)\oplus H^{0,2}(X_t,\bC)
$$
and the induced decomposition of the Gauss-Manin connection acting on $H^2$
$$
\nabla=\pmatrix{\nabla^{2,0}&*&*\cr *&\nabla^{1,1}&*\cr *&*&\nabla^{0,2}\cr}.
$$
Here the stars indicate suitable bundle morphisms -- actually with the 
lower left and upper right starts being zero by Griffiths' transversality 
property, but we do not care here. The notation $\nabla^{p,q}$ stands for
the induced (real analytic, not necessarily flat) connection on the
subbundle $t\mapsto H^{p,q}(X_t,\bC)$. The main result of this section
is the following observation.
\medskip

\noindent
{\bf 5.1. Theorem.} {\sl Let $\cX\to S$ be a deformation of compact K\"ahler
manifolds over an irreducible base~$S$. Then there exists a countable union 
$S'=\bigcup S_\nu$ of analytic subsets $S_\nu\subsetneq S$, such that 
the K\"ahler cones
$\cK_t\subset H^{1,1}(X_t,\bC)$ are invariant over $S\ssm S'$ under 
parallel transport with respect to the $(1,1)$-component $\nabla^{1,1}$ of
the Gauss-Manin connection.}
\medskip

Of course, one can apply again the theorem on each stratum $S_\nu$ instead 
of $S$ to see that there is a countable stratification of $S$ such
that the K\"ahler cone is essentially ``independent of $t$'' on each
stratum. Moreover, we have semi-continuity in the sense that $\cK_{t_0}$,
$t_0\in S'$, is always contained in the limit of the nearby cones
$\cK_t$, $t\in S\ssm S'$.
\medskip

\noindent
{\sl Proof.} The result is local over $S$, so we can possibly shrink $S$
to avoid any global monodromy (i.e., we assume that the locally
constant systems $R^k\pi_\star(\bZ_\cX)$ are constant). We then define 
the $S_\nu$'s to be the images in $S$ of those connected components of 
$C^p(\cX/S)$ which do not project onto $S$. By the fact that the projection 
is proper on each component, we infer that $S_\nu$ is an analytic subset 
of $S$. The definition of the $S_\nu$'s imply that the cohomology classes
induced by the analytic cycles $\{[Z]\}$, $Z\subset X_t$, remain exactly 
the same for all $t\in S\ssm S'$.

Since $S$ is irreducible and $S'$ is a countable union of analytic sets, 
it follows that $S\ssm S'$ is arcwise connected by piecewise smooth analytic 
arcs. Let 
$$
\gamma:[0,1]\to S\ssm S',\qquad u\mapsto t=\gamma(u)
$$ 
be such a smooth arc,
and let $\alpha(u)\in H^{1,1}(X_{\gamma(u)},\bR)$ be a family of real 
$(1,1)$-cohomology classes which are constant by parallel transport
under $\nabla^{1,1}$ (any such family is obtained by fixing $\alpha(0)$, say,
and solving the ordinary differential equation $\nabla^{1,1}\alpha=0$ on the 
interval $[0,1]$). This is equivalent to assuming that 
$$
\nabla\alpha(u)\in H^{2,0}(X_{\gamma(u)},\bC)\oplus 
H^{0,2}(X_{\gamma(u)},\bC)
$$
for all $u$. Suppose that $\alpha(0)$ is a numerically positive class
in~$X_{\gamma(0)}$. We then have
$$
\alpha(0)^p\cdot \{[Z]\}=\int_Z \alpha(0)^p>0
$$
for all $p$-dimensional analytic cycles $Z$ in $X_{\gamma(0)}$. Let us 
denote by 
$$
\zeta_Z(t)\in H^{2q}(X_t,\bZ),\qquad q=\dim X_t -p,
$$ 
the family of cohomology classes equal to $\{[Z]\}$ at $t=\gamma(0)$,
such that $\nabla\zeta_Z(t)=0$ (i.e.\ constant with respect to the Gauss-Manin
connection).  By the above discussion, $\zeta_Z(t)$ is of type
$(q,q)$ for all $t\in S$, and when $Z\subset X_{\gamma(0)}$ varies,
$\zeta_Z(t)$ generates all classes of analytic cycles in $X_t$ if
$t\in S\ssm S'$.  Since $\zeta_Z$ is $\nabla$-parallel and
$\nabla\alpha(u)$ has no component of type $(1,1)$, we find
$$
{d\over du}(\alpha(u)^p\cdot\zeta_Z(\gamma(u)) =
p \alpha(u)^{p-1}\cdot \nabla\alpha(u)\cdot \zeta_Z(\gamma(u))= 0.
$$
We infer from this that $\alpha(u)$ is a numerically positive class
for all $u\in[0,1]$. This argument shows that the set $\cP_t$ of 
numerically positive classes in $H^{1,1}(X_t,\bR)$ is invariant
by parallel transport under $\nabla^{1,1}$ over $S\ssm S'$.

By a standard result of Kodaira-Spencer [KS60] relying on elliptic PDE
theory, every K\"ahler class in $X_{t_0}$ can be deformed to a nearby 
K\"ahler class in nearby fibres~$X_t$. This implies that the set of 
$t\in S\ssm S'$ for which a given connected component of $\cP_t^\circ$ 
coincides with the K\"ahler cone $\cK_t$ is open. As $S\ssm S'$ is connected,
these sets must be either empty or equal to $S\ssm S'$, hence the connected
component of $\cP_t^\circ$ which yields the K\"ahler cone remains the same
for all $t\in S\ssm S'$. The theorem is thus proved (notice moreover that the 
remark concerning the semi-continuity of K\"ahler cones stated after
Theorem 5.1 follows from the result by Kodaira-Spencer).\qed \medskip

From the above results, one can hope for a much stronger semi-continuity 
statement than the one stated by Kodaira-Spencer. Namely,
we make the following conjecture, which we will consider in a 
forthcoming paper.

\medskip

\noindent{\bf 5.2. Conjecture.} {\sl Let $\cX\to S$ be a deformation of
compact complex manifolds over an irreducible base $S$. Assume that
one of the fibres $X_{t_0}$ is K\"ahler. Then there exists a countable
union $S'\subsetneq S$ of analytic subsets in the base such that
$X_t$ is K\"ahler for $t\in S\ssm S'$. Moreover, $S'$ can be chosen so that
the K\"ahler cone is invariant over $S\ssm S'$, under parallel transport 
by the Gauss-Manin connection.}
\medskip

In other words, the K\"ahler property should be open for the countable
Zariski topology on the base $S$ (as well as for the usual metric space 
topology on $S$, by Kodaira-Spencer). We are not sure what to think about the 
remaining fibres $X_t$, $t\in S'$, but a natural expectation would be that
they are in the Fujiki class $\cC$, at least under the assumption that Hodge
decomposition remains valid on those fibres -- which is anyway 
a necessary condition for the expectation to hold true.
\vfill\eject

\noindent {\twelvebf References}
\bigskip

{\eightpoint
\parindent = 1.5cm

\item{\bf [Bar75]} Barlet, D.\ --- {\sl Espace analytique r\'eduit des 
cycles analytiques complexes compacts d'un espace analytique complexe 
de dimension finie}, Fonctions de plusieurs variables complexes, II 
(S\'em.\ François Norguet, 1974--1975), Lecture Notes in Math., 
Vol.~482, Springer, Berlin (1975) 1--158. 

\item{\bf [Bon93]} Bonavero, L.\ --- {\sl In\'egalit\'es de Morse
holomorphes singuli\`eres}, C.\ R.\ Acad.\ Sci.\ S\'erie I {\bf 317} (1993),
1163--1166.

\item{\bf [Bon98]} Bonavero, L.\ --- {\sl In\'egalit\'es de morse holomorphes 
singuli\`eres}, J.\ Geom.\ Anal.\ {\bf 8} (1998), 409--425.

\item{\bf [Bou89]} Bouche, T.\ --- {\sl In\'egalit\'es de Morse 
pour la $d^{\prime \prime}$-cohomologie sur une vari\'et\'e
non-compacte}, Ann.\ Sci.\ Ecole Norm.\ Sup.\  {\bf 22} 
(1989), 501--513.

\item{\bf [Buc99]} Buchdahl, N.\ --- {\sl On compact K\"ahler surfaces}, 
Ann.\ Inst.\ Fourier {\bf 49} (1999) 287--302.

\item{\bf [Buc00]} Buchdahl, N.\ --- {\sl A Nakai-Moishezon criterion for 
non-K\"ahler surfaces}, Ann.\ Inst.\ Fourier {\bf 50} (2000) 1533--1538.

\item{\bf [CP90]} Campana, F., Peternell, T.\ --- {\sl Algebraicity of the 
ample cone of projective varieties} J.\ Reine Angew.\ Math.\ {\bf 407} (1990).

\item{\bf [Dem85a]} Demailly, J.-P.\ --- {\sl Mesures de Monge-Amp\`ere
et caract\'erisation g\'eom\'etrique des vari\'et\'es alg\'e\-briques
affines}, M\'em.\ Soc.\ Math.\ France (N.S.) {\bf 19} (1985) 1--124.

\item{\bf [Dem85b]} Demailly, J.-P.\ --- {\sl Champs magn\'etiques et
in\'egalit\'es de Morse pour la $d^{\prime \prime}$-cohomo\-logie},
Ann.\ Inst.\ Fourier\ {\bf 35} (1985).

\item{\bf [Dem90a]} Demailly, J.-P.\ --- {\sl Cohomology of $q$-convex 
Spaces in Top Degrees}, Math.\ Z.\ {\bf 204} (1990), 283--295.

\item{\bf [Dem90b]} Demailly, J.-P.\ --- {\sl Singular hermitian metrics on 
positive line bundles}, Proceedings of the Bayreuth conference ``Complex 
algebraic varieties'', April~2-6, 1990, edited by K.~Hulek, T.~Peternell, 
M.~Schneider, F.~Schreyer, Lecture Notes in Math.\ ${\rm n}^\circ\,$1507, 
Springer-Verlag, 1992.

\item{\bf [Dem92]} Demailly, J.-P.\ --- {\sl Regularization of closed
positive currents and intersection theory}, J.\ Algebraic Geometry,
{\bf 1} (1992), 361--409.

\item{\bf [Dem93]} Demailly, J.-P.\ --- {\sl A numerical criterion for
very ample line bundles}, J.\ Diff.\ Geom.\ {\bf 37} (1993), 323--374.

\item{\bf [DEL99]} Demailly, J.-P., Ein, L., Lazarsfeld, R.\ --- {\sl 
A subadditivity property of multiplier ideals}, math.AG/0002035,
Michigan Math.\ J., special volume in honor of William Fulton,
{\bf 48} (2000), 137--156.

\item{\bf [DPS94]} Demailly, J.-P., Peternell, Th., Schneider, M.\ --- {\sl
Compact complex manifolds with numerically effective tangent bundles},
J.\ Algebraic Geometry {\bf 3} (1994) 295--345.

\item{\bf [DPS00]} Demailly, J.-P., Peternell, Th., Schneider, M.\ --- {\sl
Pseudo-effective line bundles on compact K\"ahler manifolds}, to appear
in the International Journal of Mathematics.

\item{\bf [Eys00]} Eyssidieux, P.\ --- {\sl Th\'eor\`emes de 
Nakai-Moishezon pour certaines classes de type $(1,1)$ et applications}, 
Pr\'epublication Universit\'e de Toulouse III, 2000.

\item{\bf [Huy99]} Huybrechts, D.\ --- {\sl Compact Hyperk\"ahler Manifolds:
Basic Results}, Invent.\ math.\ {\bf 135} (1999), 63--113.

\item{\bf [Huy01]} Huybrechts, D.\ --- {\sl The Projectivity Criterion 
for Hyperk\"ahler manifolds as a Consequence of the Demailly-Paun Theorem},
personal communication, manuscript K\"oln University, May 2001, 3~p.

\item{\bf [JS93]} Ji, S.,  Shiffman, B.\ --- {\sl Properties of compact 
complex manifolds carrying closed positive currents}, J.\ Geom.\ Anal.\ 
{\bf 3}, (1993) 37--61.

\item{\bf [Kle66]} Kleiman, S.\ --- {\sl Toward a numerical theory of
ampleness}, Ann.\ Math.\ {\bf 84} (1966), 293--344. 
Kodaira, K.; Spencer, D. C. 

\item{\bf [KS60]} Kodaira, K., Spencer, D.C.\ {\sl On deformations of 
complex analytic structures. III. Stability theorems for complex structures},
Ann.\ of Math.\ {\bf 71} (1960) 43--76. 

\item{\bf [Lam99a]} Lamari, A.\ --- {\sl Courants k\"ahl\'eriens et 
surfaces compactes}, Ann.\ Inst.\ Fourier {\bf 49} (1999) 263--285.

\item{\bf [Lam99b]} Lamari, A.\ --- {\sl Le cone K\"ahl\'erien d'une 
surface}, J.\ Math.\ Pures Appl.\ {\bf 78} (1999) 249--263

\item{\bf [Pau98a]} Paun, M.\ --- {\sl Sur l'effectivit\'e num\'erique 
des images inverses de fibr\'es en droites}, Math.\ Ann.\ {\bf 310} (1998),
411--421.

\item{\bf [Pau98b]} Paun, M.\ --- {\sl Fibr\'e en droites num\'eriquement 
effectifs et vari\'et\'es K\"ahl\'eriennes compactes 
\`a courbure de Ricci nef}, Th\`ese, Universit\'e de Grenoble 1 (1998), 
80$\,$p.

\item{\bf [Pau00]} Paun, M.\ --- {\sl Semipositive (1,1)--cohomology 
classes on projective manifolds}, Preprint de l'I.R.M.A.

\item{\bf [Siu74]} Siu, Y.-T.\ --- {\sl Analyticity of sets associated to
Lelong numbers and the extension of closed positive currents},  
Invent.\ Math.\ {\bf 27} (1974).

\item{\bf [Siu84]} Siu, Y.-T.\ --- {\sl Some recent results in complex
manifold theory for the semi-positive case}, survey article in the 
Proceedings of the international congress held in Bonn, 1984. 

\item{\bf [Yau78]} Yau, S.-T.\ --- {\sl On the Ricci curvature of a 
complex K\"ahler manifold and the complex Monge--Amp\`ere equation},
Comm.\ Pure Appl.\ Math.\ {\bf 31} (1978). 

\vskip30pt
\noindent
Jean-Pierre Demailly\\
Universit\'e de Grenoble I, Institut Fourier, UMR 5582 du CNRS\\
BP74, 100 rue des Maths, 38402 Saint-Martin d'H\`eres Cedex, France\\
{\sl E-mail:}\/ demailly@fourier.ujf-grenoble.fr
\medskip

\noindent
Mihai Paun\\
Universit\'e Louis Pasteur, D\'epartement de Math\'ematiques\\
7, rue Ren\'e Descartes, 67084 Strasbourg Cedex, France\\
{\sl E-mail:}\/ paun@math.u-strasbg.fr
\medskip

\noindent
(September 14, 2001; printed on \today)

}

\end